\documentclass[1p]{elsarticle}
\usepackage[final]{epsfig}
\usepackage{graphics}
\usepackage{amsmath}
\usepackage{amsthm}
\usepackage{amsfonts}
\usepackage{latexsym}
\usepackage{amssymb}
\usepackage{graphicx}
\usepackage{color}
\newtheorem{lemma}{Lemma}[section]
\newtheorem{proposition}[lemma]{Proposition}
\newtheorem{remark}[lemma]{Remark}
\newtheorem{example}[lemma]{Example}
\newtheorem{theorem}[lemma]{Theorem}
\newtheorem{definition}[lemma]{Definition}

\begin{document}
\newcommand{\eps}{{\varepsilon}}
\newcommand{\proofend}{$\Box$\bigskip}
\newcommand{\C}{{\mathbf C}}
\newcommand{\Q}{{\mathbf Q}}
\newcommand{\R}{{\mathbf R}}
\newcommand{\Z}{{\mathbf Z}}
\newcommand{\RP}{{\mathbf {RP}}}
\def\proof{\paragraph{Proof.}}

\def\Re{{\mbox{Re}}}
\def\Im{{\mbox{Im}}}
\def\dfrac#1#2{\displaystyle{\frac{#1}{\,#2\,}}}
\def\L{{L}}
\def\A{{A}}
\def\div{{\operatorname{div}}}

\def\Dx4{{\left.D_{x_{\phantom{j}}}\!\!\right|_{w=4}}}
\def\wLap{{L_{_\Delta}}}
\def\Lap{{\left.\wLap_{_{\phantom{j}}}\!\!\right|_{w=4}}}
\def\Lapf{{\left.\wLap_{_{\phantom{j}}}\!\!\right|_{w=4}}}
\def\Lapt{{\left.\wLap_{_{\phantom{j}}}\!\!\right|_{w=2}}}
\def\Lapi{{\left.\wLap_{_{\phantom{j}}}\!\!\right|_{w=\infty}}}
\def\Lapz{{\left.\wLap_{_{\phantom{j}}}\!\!\right|_{w=0}}}
\newcommand{\dir}{\mathbf{t}}
\newcommand{\tmax}{\mathbf{t}_{\rm max}}
\newcommand{\tmin}{\mathbf{t}_{\rm min}}
\newcommand{\kmax}{k_{\rm max}}
\newcommand{\kmin}{k_{\rm min}}
\newcommand{\n}{\mathbf{n}}
\newcommand{\p}{\mathbf{p}}
\newcommand{\q}{\mathbf{q}}
\newcommand{\surf}{\mathcal{S}}

\def\x{\textit{\textbf{x}}}
\def\y{\textit{\textbf{y}}}
\def\e{\textit{\textbf{e}}}

\newcounter{ab}
\newcommand{\ab}[1]
{\stepcounter{ab}$^{\bf AB\theab}$%
\footnotetext{\hspace{-3.7mm}$^{\blacksquare\!\blacksquare}$
{\bf AB\theab:~}#1}}

\newcounter{bk}
\newcommand{\bk}[1]
{\stepcounter{bk}$^{\bf BK\thebk}$%
\footnotetext{\hspace{-3.7mm}$^{\blacksquare\!\blacksquare}$
{\bf BK\thebk:~}#1}}

\newcounter{st}
\newcommand{\st}[1]
{\stepcounter{st}$^{\bf ST\thest}$%
\footnotetext{\hspace{-3.7mm}$^{\blacksquare\!\blacksquare}$
{\bf ST\thest:~}#1}}


\begin{frontmatter}

\title{Discrete spherical means of directional derivatives \\ and Veronese maps}
\author[label1]{Alexander Belyaev\,}
\author[label2]{Boris Khesin\,}
\author[label3]{Serge Tabachnikov\,}
\address[label1]{Electrical, Electronic \& Computer Engineering,
Heriot-Watt University, Edinburgh, EH14 4AS, UK}
\address[label2]{Department of Mathematics, University of Toronto,
 Toronto, ON M5S 2E4, Canada\,}
\address[label3]{Department of Mathematics, Pennsylvania State University,
 University Park, PA 16801, USA\,}

\date{}

\begin{abstract}
We describe and study geometric properties of discrete circular and spherical
means of directional derivatives of functions, as well as discrete approximations
of higher order differential operators.
For an arbitrary dimension we present a general
construction for obtaining discrete spherical means of directional
derivatives. The construction is based on using the Minkowski's existence
theorem and Veronese maps. Approximating the directional derivatives
by appropriate finite differences allows one to obtain finite difference
operators with good rotation invariance properties.
In particular, we use discrete circular and spherical
means to derive  discrete approximations of various linear and nonlinear first-{}
and second-order differential operators, including discrete Laplacians.
A practical potential of our
approach is demonstrated by considering applications to nonlinear
filtering of digital images and surface curvature estimation.
\end{abstract}

\end{frontmatter}
\section*{Introduction}

Mean value properties of functions play crucial role in analysis
of many partial differential equations \cite{Friedman-Littman_62},
numerical interpolation and integration \cite{Stroud_book71},
computer tomography \cite{Natterer_siam01}, and many other areas
of mathematics and engineering. In this paper we undertake a study
of discrete circular and spherical means of directional derivatives
of functions in view of their applications to the finite  difference
methods. We present a general construction for obtaining discrete
spherical means of directional derivatives, analyze general properties
of such means, and use them to derive finite difference approximations
with good rotation-invariant properties for linear and nonlinear
first-{} and second-order partial differential operators.
The two principal components of our approach are based
on the use of the Minkowski's existence theorem
\cite{Schneider_93} and Veronese maps \cite{Harris_ag}.

For the sake of simplicity, we focus on finite difference approximations
of very basic differential operators, the gradient and Laplacian. We
start from an elementary analysis of discrete circular means of derivatives
and demonstrate how such discrete means can be used for practical
estimation of surface curvatures (Section \ref{sec-a}). Then, in any
dimension, we obtain Minkowski-type formulas for general and regular
grids by using the Veronese map (Section \ref{sec-t}).
Finally, we analyze
rotation invariance properties of discrete 2D gradient, Laplacian, and
related nonlinear operators defined on a square grid (Section~\ref{sec-l}).

Various finite difference approximations of the gradient and Laplacian
are widely used in fluid mechanics, electromagnetics, finance modeling,
image processing, and many other areas. Our interest in constructing
reliable discrete approximations of these operators is threefold.
Firstly, reliable and consistent approximations of the gradient
and Laplacian are needed in various diffusion-type hydrodynamical models,
such as the convection-diffusion equation
$$
u_t+v\cdot\nabla u=\epsilon\Delta u\,,
$$
where the density $u$ is convected by a velocity field $v$ and dissipates at different time scales. Its attractors exhibit very peculiar behavior  as
$t\to \infty$ and $\epsilon\to 0$  (see \cite{Arnold-Khesin98})
and we hope that the approximation scheme developed in this paper could complement the analytic treatment of the attractors.
The Laplacian case is also particularly useful in computations related
to the fast dynamo problem, where the proposed approximation formulas may simplify
computations of the magnetic diffusion for iterations of a seed magnetic
field, see \cite{Childress}.

Secondly, linear and nonlinear
diffusion-type equations are widely used in modern image processing
for a variety of tasks including edge detection, image restoration,
and image decomposition into structure and and texture components
\cite{Weickert_book98,Aubert-Kornprobst_mpip2e06}.
As demonstrated in \cite{Weickert_jvcir02}, the use of discrete
nonlinear diffusion operators with good rotation invariance properties
is highly beneficial for such tasks. In addition, nonlinear image diffusion provides us with a useful testbed for the theoretical considerations below.
Thirdly, accurate and reliable approximations of the surface gradient
and Laplacian are key ingredients in many shape analysis studies. Below we discuss applications of the approximation formulas both in image processing and geometry.
Finally we note that the approach via the Veronese map described in the present paper
can also be used for higher order differential operators, such as bi-Laplacian.


\section{Circular means of derivatives and their applications}\label{sec-a}

\paragraph*{\bf Continuous and discrete circular means}
A simple way to demonstrate the rotational invariance of
the 2D Laplacian
$$
\Delta\equiv\frac{\partial^2}{\partial x^2}+\frac{\partial^2}{\partial y^2}
$$
consists of representing it as the circular
mean (respectively, the spherical mean in 3D) of the second-order directional derivatives
\begin{equation}\label{mv-lap}
\frac12\Delta f=\dfrac{1}{\pi}\int_{0}^{\pi}
\dfrac{\partial^2 f}{\partial\e_\varphi^2}d\varphi,
\quad\mbox{where }
\e_\varphi=\left(\cos\varphi,\sin\varphi\right),\phantom{l}
\dfrac{\partial}{\partial {\e_\varphi}}=
\cos\varphi\,\dfrac{\partial}{\partial x}+\sin\varphi\,\dfrac{\partial}{\partial y}.
\end{equation}
The rotational invariance of the squared gradient $|\nabla f|^2$ can be
demonstrated in a similar way as
\begin{equation}\label{mv-grad2}
\frac12|\nabla f|^2=\dfrac{1}{\pi}\int_{0}^{\pi}
\left(\dfrac{\partial f}{\partial {\e_\varphi}}\right)^2d\varphi.
\end{equation}
More generally, the mean value representation for a directional derivative of $f$
can be obtained from
\begin{equation}\label{mv-grad1}
\frac12\left[a\dfrac{\partial f}{\partial x} + b\dfrac{\partial f}{\partial y}\right]
=\dfrac{1}{\pi}\int_{0}^{\pi}
\left[a\cos\varphi+b\sin\varphi\right]\dfrac{\partial f}{\partial {\e_\varphi}}\,d\varphi\,,
\end{equation}
where $a$ and $b$ are constants.
Notice that (\ref{mv-grad1}) with $a=\partial/\partial x$ and
$b=\partial/\partial y$ implies (\ref{mv-lap}) and setting $a=\partial f/\partial x$
and $b=\partial f/\partial y$ yields the equality (\ref{mv-grad2}).

Let us now derive discrete counterparts of the above mean value representations
(\ref{mv-lap}), (\ref{mv-grad2}), and (\ref{mv-grad1}). Consider a set of directions
$\e_k=(\cos\varphi_k,\sin\varphi_k)$, $k=1,2,\dots,n$. Then a weighted sum of the second-order
directional derivatives can be expressed as
\begin{equation}\label{sum}
\begin{array}{l}
\sum w_k\dfrac{\partial^2f}{\partial\e^2_{\varphi_k}}
=\dfrac{\partial^2 f}{\partial x^2}\sum w_k\cos^2\varphi_k+
\dfrac{\partial^2 f}{\partial y^2}\sum w_k\sin^2\varphi_k
\vspace*{2pt} \phantom{MMM}\\ \phantom{MMM}
{}+2\dfrac{\partial^2 f}{\partial x\partial y}\sum w_k\cos\varphi_k\sin\varphi_k
=\dfrac12\Delta f\sum w_k
\vspace*{6pt} \\ \phantom{MMM}
{}+\dfrac12\left(\dfrac{\partial^2 f}{\partial x^2}-\dfrac{\partial^2 f}{\partial y^2}\right)
\Re\left(\sum w_ke^{2i\varphi_k}\right)+
\dfrac{\partial^2 f}{\partial x\partial y}
\Im\left(\sum w_ke^{2i\varphi_k}\right).
\end{array}
\end{equation}
Finding weights $\{w_k\}$ such that
\begin{equation}\label{wk}
\sum w_ke^{2i\varphi_k}=0
\end{equation}
yields, after a simple normalization
$
\{w_k\}\longrightarrow\left\{w_k\left/\,\sum w_k\right.\right\},
$
a discrete mean value representation for $\Delta f$:
\begin{equation}\label{dmv4lap}
\dfrac12\Delta f = \sum w_k\dfrac{\partial^2 f}{\partial\e_{\varphi_k}^2},
\end{equation}
and similar representations
of the squared gradient
\begin{equation}\label{dmv4grad}
\dfrac12|\nabla f|^2 = \sum w_k\left(\dfrac{\partial f}{\partial \e_{\varphi_k}}\right)^2
\end{equation}
and directional derivative
\begin{equation}\label{dmv4dirder}
\dfrac12\left[a\dfrac{\partial f}{\partial x}+b\dfrac{\partial f}{\partial y}\right]
=\sum w_k\left[a\cos\varphi_k+b\sin\varphi_k\right]
\dfrac{\partial f}{\partial\e_{\varphi_k}},
\end{equation}
respectively.
Obviously the same set of weights $\{w_k\}$ can be used to obtain
a discrete mean value representation for a quasi-Laplacian
\begin{equation}\label{qLap}
\dfrac12\nabla\cdot(a(x,y)\nabla f)=\sum w_k\dfrac{\partial}{\partial \e_{\varphi_k}}
\left(a(x,y)\dfrac{\partial f}{\partial\e_{\varphi_k}}\right).
\end{equation}
Note that the appearance of double angles in the formula (\ref{wk}) exactly
corresponds to the Veronese map $V:(\cos \phi, \sin \phi) \mapsto (\cos 2\phi, \sin 2\phi)$,
considered in Example~\ref{example2d} of the next section. 

\paragraph*{\bf Discrete mean value formulas and Minkowski problem for polygons}

The key problem of finding normalized weights $\{w_k~|~\sum w_k=1\}$ in the relation $\sum w_ke^{2i\varphi_k}=0$
(see (\ref{wk}) above) has a simple geometric meaning
and is directly related to the Minkowski existence theorem (also called
Minkowski's Problem, see e.g. \cite{Schneider_93}) for polyhedra.
\begin{theorem}[Minkowski's existence theorem] \label{Mink-exist}
Consider unit vectors $\n_1,\n_2,\dots,\n_k$ spanning $\R^m$
and a set of positive weights $w_1,w_2,\dots,w_k$. Then the equation
$$
w_1\n_1+w_2\n_2+\dots+w_k\n_k=0
$$
is a necessary and sufficient condition for existence of a convex polyhedron
with facet unit normals $\n_1,\n_2,\dots,\n_k$ and corresponding facet areas
$w_1,w_2,\dots,w_k$. Furthermore, this polyhedron is unique up to translation.
\end{theorem}
The necessity part of this theorem is, of course, trivial and follows
immediately from the  Gauss divergence theorem.

\begin{figure}[htbp]
\centering
\includegraphics[width=0.3\textwidth]{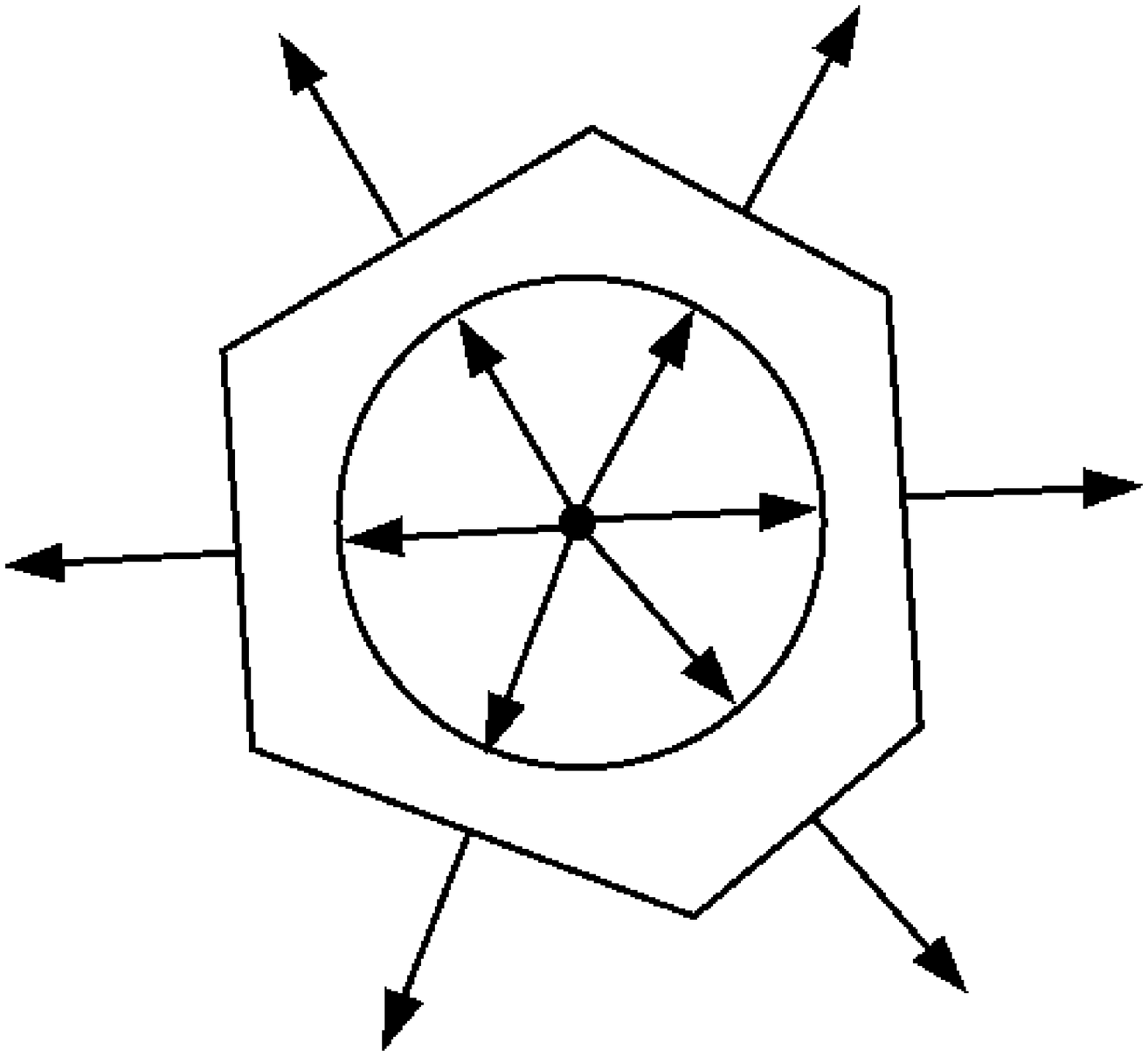}
\quad
\includegraphics[width=0.2\textwidth]{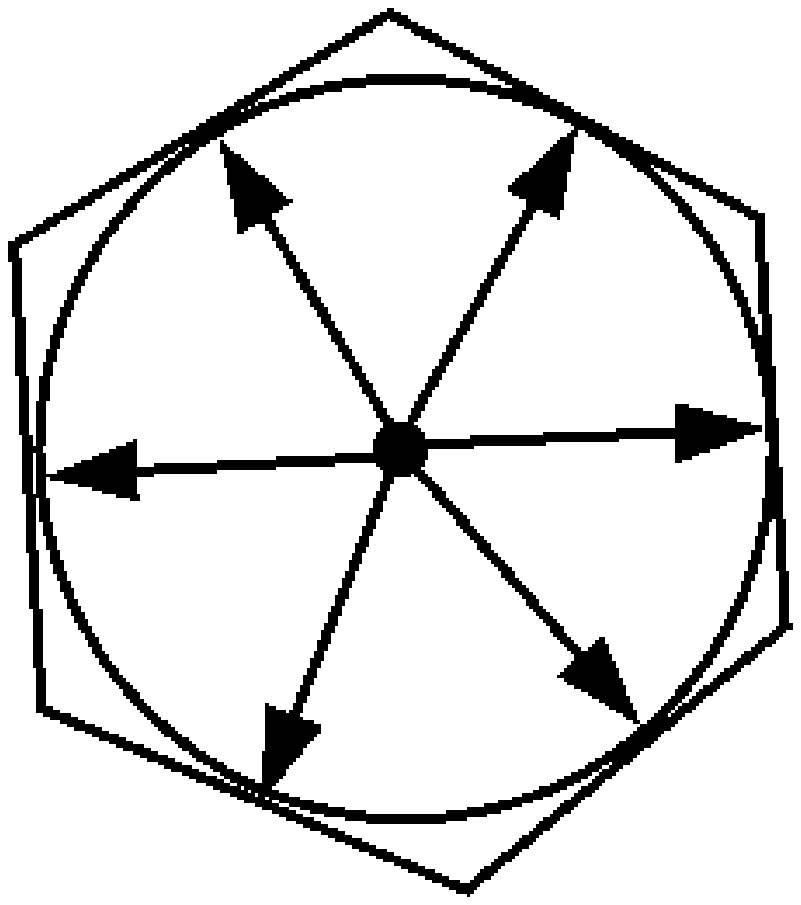}
\qquad
\includegraphics[width=0.25\textwidth]{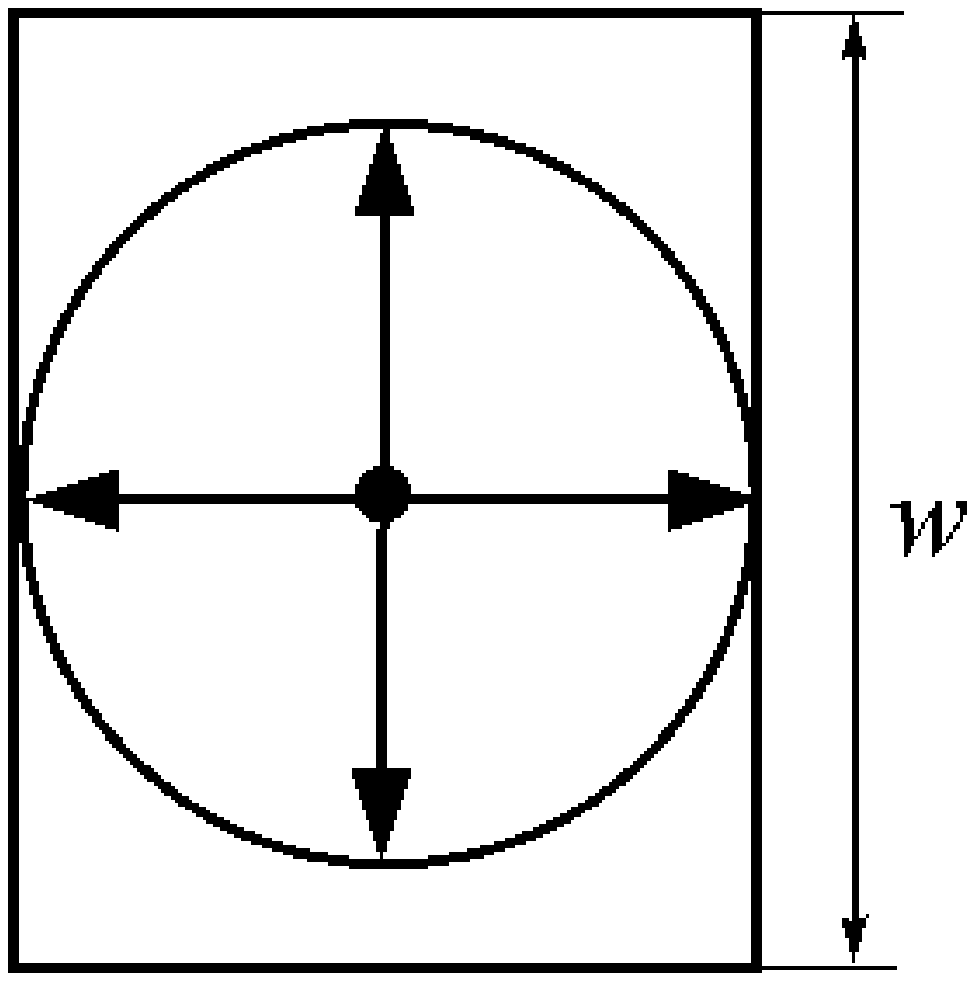}
\caption{Left: Minkowski's existence theorem for polygons provides
us with a geometric interpretation of (\ref{wk}). Middle: a special case
of a polygon circumscribed around the unit circle was considered
in \cite{Langer_cagd07} in connection with the curvature estimation
approach proposed there. Right: $3\times3$ stencils
for estimating the gradient and Laplacian of functions defined
on a square grid contain four directions and correspond to rectangles.
}
\label{Mink}
\end{figure}

Let us note that in (\ref{mv-lap}), (\ref{mv-grad2}), and (\ref{mv-grad1})
one deals with
the directions $\e(\varphi)=(\cos\varphi,\sin\varphi)$, $0\leq\varphi<\pi$,
which can be parameterized by the unit vector $e^{2i\varphi}$.
Now consider a closed polygon such that the vectors $e^{2i\varphi_k}$
are the outward normals of the polygon's edges. Then, according to the Gauss divergence theorem, the edge lengths provide us with the desired set of weights $\{w_k\}$. Vice versa, given
a set of unit vectors $\{e^{2i\varphi_k}\}$ and positive weights
$\{w_k\}$ satisfying the first condition in (\ref{wk}), Minkowski's existense theorem guarantees an existence of a convex polygon
with edge lengths $\{w_k\}$ and outward normals $e^{2i\varphi_k}$.
The left image of Figure~\ref{Mink} illustrates this geometric interpretation
of (\ref{wk}).

In particular, if we consider a polygon circumscribed around the unit circle,
as seen in the middle image of Figure~\ref{Mink},
the lengths of the polygon sides are given by
$$
\tan\beta_{j-1}+\tan\beta_j,\quad \beta_k=\varphi_{k+1}-\varphi_k,
$$
and, therefore, the normalized weights are
\begin{equation}\label{w-circ}
w_j=\dfrac{\tan\beta_{j-1}+\tan\beta_j}{\sum\left(\tan\beta_{k-1}+\tan\beta_k\right)}.
\end{equation}
This set of weights was used in \cite{Langer_cagd07} for estimating the mean curvature
on surfaces approximated by dense triangle meshes.

In Section~\ref{sec-t} we show how to solve this key problem in any dimension and, in particular, get a discretization of the corresponding Laplacian.

\paragraph*{\bf Application to curvature estimation}
In the previous examples, we have considered function data defined
on a regular grid. Here we deal with less organized data, namely
with triangle meshes approximating smooth surfaces
{embedded into the Euclidean space $\R^{3}$}.

Reliable estimation of curvature characteristics of polygonal surfaces
is very important for many computer graphics and geometric modeling applications.
Starting from a seminal paper of Taubin \cite{Taubin_iccv95} integral-based approaches
are frequently used for curvature estimation purposes
\cite{Langer_cagd07,pottmann-2009-iir,Meigot_spm09} (see also references therein).
Our approach to curvature estimation is based on using discrete circular
means of directional surface derivatives.
Given a smooth surface approximated by a dense triangle mesh,
we employ such discrete circular means
for estimating the mean curvature $H=(\kmax+\kmin)/2$ and the so-called curvedness
\cite{Koenderink-VanDoorn_92}
$R=\sqrt{\left(\kmax^2+\kmin^2\right)/2}$,
where $\kmax$ and $\kmin$ stand for the principal curvatures.
(While the curvedness $R$ is not a classical curvature measure, it is widely used
in numerous applications including protein interaction analysis \cite{Bradford-Westhead_bio05}, heart physiology \cite{Zhong_09},
and human visual perception \cite{ValentiICCV2009}.)

Let $\n$ be a surface orientation normal. Denote by $\tmax$ and $\tmin$
the principal directions corresponding to the principal curvatures $\kmax$ and $\kmin$, respectively. For a point $P$ on the surface, consider a unit tangent vector
$\dir(\varphi)$, which makes angle $\varphi$ with $\tmax$, and the corresponding
directional curvature $k(\varphi)$.
Similar to (\ref{mv-lap}) and (\ref{mv-grad2}) we have
\begin{equation}\label{Kmean}
H\equiv\dfrac{1}{2}\left(\kmax+\kmin\right)=\dfrac{1}{2\pi}\int_0^{2\pi}
k(\varphi)\,d\varphi
\end{equation}
and
\begin{equation}\label{Ktotal}
R^2\equiv\dfrac{1}{2}\left(\kmax^2+\kmin^2\right)=\dfrac{1}{2\pi}\int_0^{2\pi}
\left(\dfrac{\partial\n}{\partial\dir(\varphi)}\right)^2\,d\varphi,
\end{equation}
respectively. Here (\ref{Kmean}) immediately follows from Euler's formula
$$
k(\varphi)=\kmax\cos\varphi^2+\kmin\sin\varphi^2
$$
and (\ref{Ktotal}) is obtained from a similar representation
$$
\left(\frac{\partial\n}{\partial\dir(\varphi)}\right)^2
=\kmax^2\cos^2\varphi+\kmin^2\sin^2\varphi,
$$
which, in  turn, can be easily derived from the Rodrigues formulas
$$
\dfrac{\partial\n}{\partial\tmax}=-\kmax\,\tmax,\quad
\dfrac{\partial\n}{\partial\tmin}=-\kmin\,\tmin.
$$

Similar to (\ref{dmv4lap}) and (\ref{dmv4grad}), let us consider
discrete counterparts of (\ref{Kmean}) and (\ref{Ktotal})
\begin{equation}\label{HR-discrete}
H=\sum w_jk(\varphi_j)\quad\mbox{and}\quad
R^2=\sum w_j\left(\dfrac{\partial\n}{\partial\dir(\varphi_j)}\right)^2,
\end{equation}
where weights $\{w_j\}$
{satisfy (\ref{wk})
and are normalized by $\sum w_k=1$.}
For the sake of simplicity, we use the set of weights (\ref{w-circ})
proposed in \cite{Langer_cagd07}.

\begin{figure}[htbp]
\centering
\includegraphics[width=0.6\textwidth]{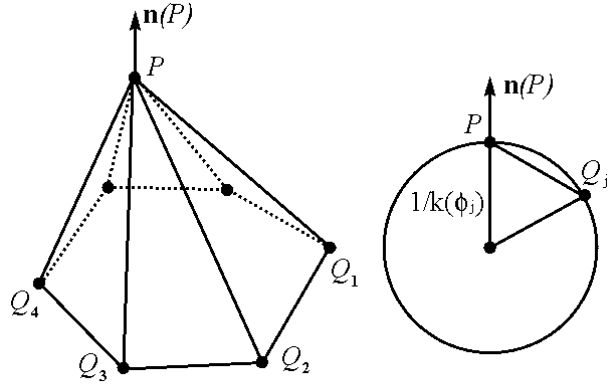}
\caption{
Left: a mesh vertex $P$ and its 1-ring neighboring vertices $\{Q_j\}$.
Right: a geometric explanation of the directional curvature approximation
in (\ref{HR-approx}).
}
\label{nghbrs}
\end{figure}

Consider now a triangle mesh approximating the surface and assume that
each mesh vertex $P$ is equipped with a unit normal $\n(P)$ which
delivers a fairly good approximation of the corresponding surface normal.
Given a mesh vertex $P$ and its 1-ring neighboring vertices $\{Q_j\}$,
as seen in the left image of Figure~\ref{nghbrs},
simple approximations of tangential vectors at $P$ and the corresponding
directional curvature $k(\varphi_j)$ and directional derivative of $\n$
are given by
\begin{equation}\label{HR-approx}
\dir_j\approx\dfrac{\vec{PQ_j}}{\|P-Q_j\|},
\quad
k(\varphi_j)\approx\dfrac{2\vec{PQ_j}\cdot\n}{\|P-Q_j\|^2},
\quad
\dfrac{\partial\n}{\partial\dir_j}\approx\frac{\n(Q_j)-\n(P)}{\|Q_j-P\|},
\end{equation}
respectively. We refer to the right image of Figure~\ref{nghbrs} for
a geometric explanation of the directional curvature approximation.
Finally, (\ref{HR-discrete}) and (\ref{HR-approx}) provide us
with discrete approximations of the mean curvatures $H$ and curvedness $R$.

\begin{figure}[htbp]
\centering
\includegraphics[width=0.32\textwidth]{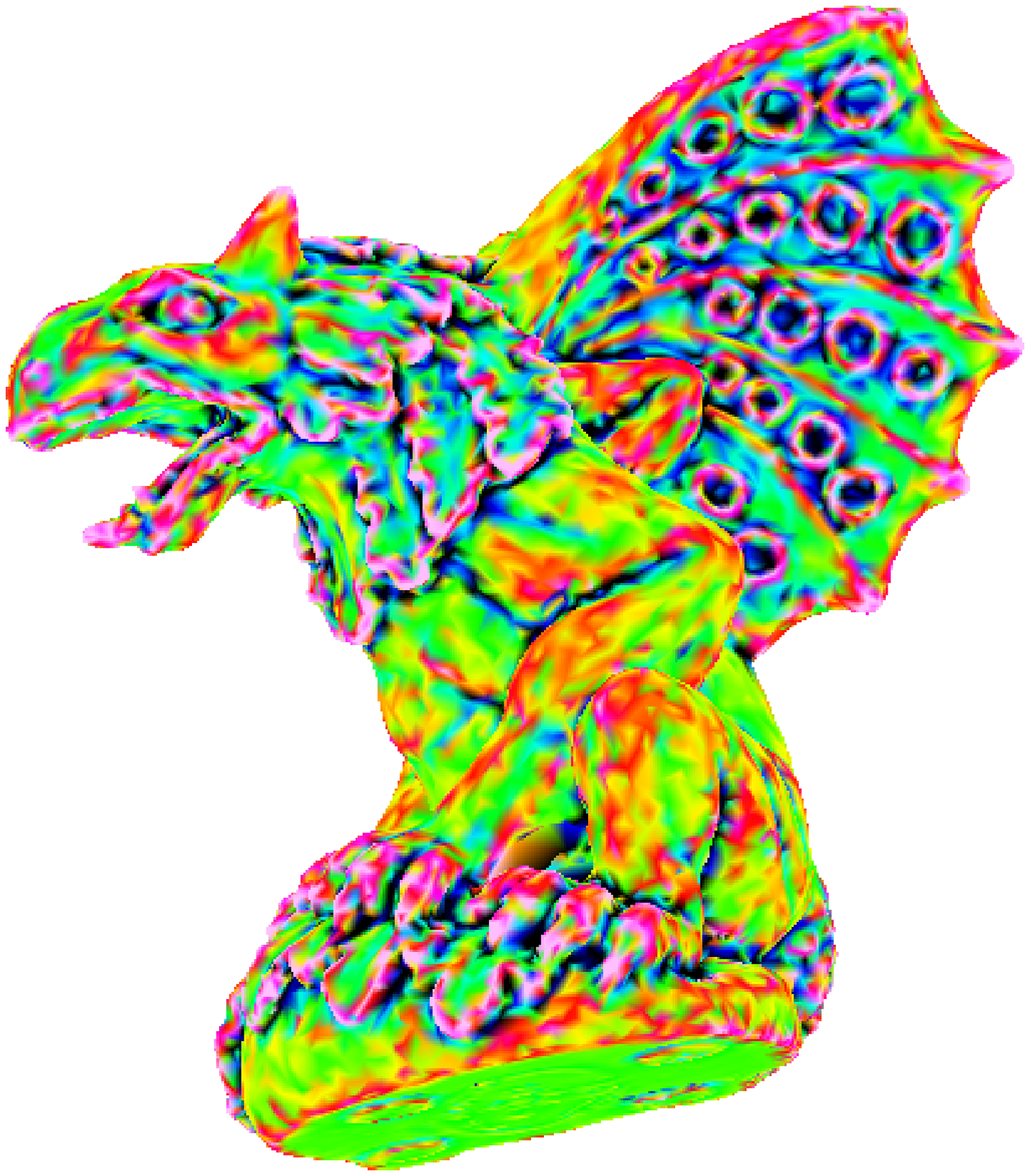}
\includegraphics[width=0.32\textwidth]{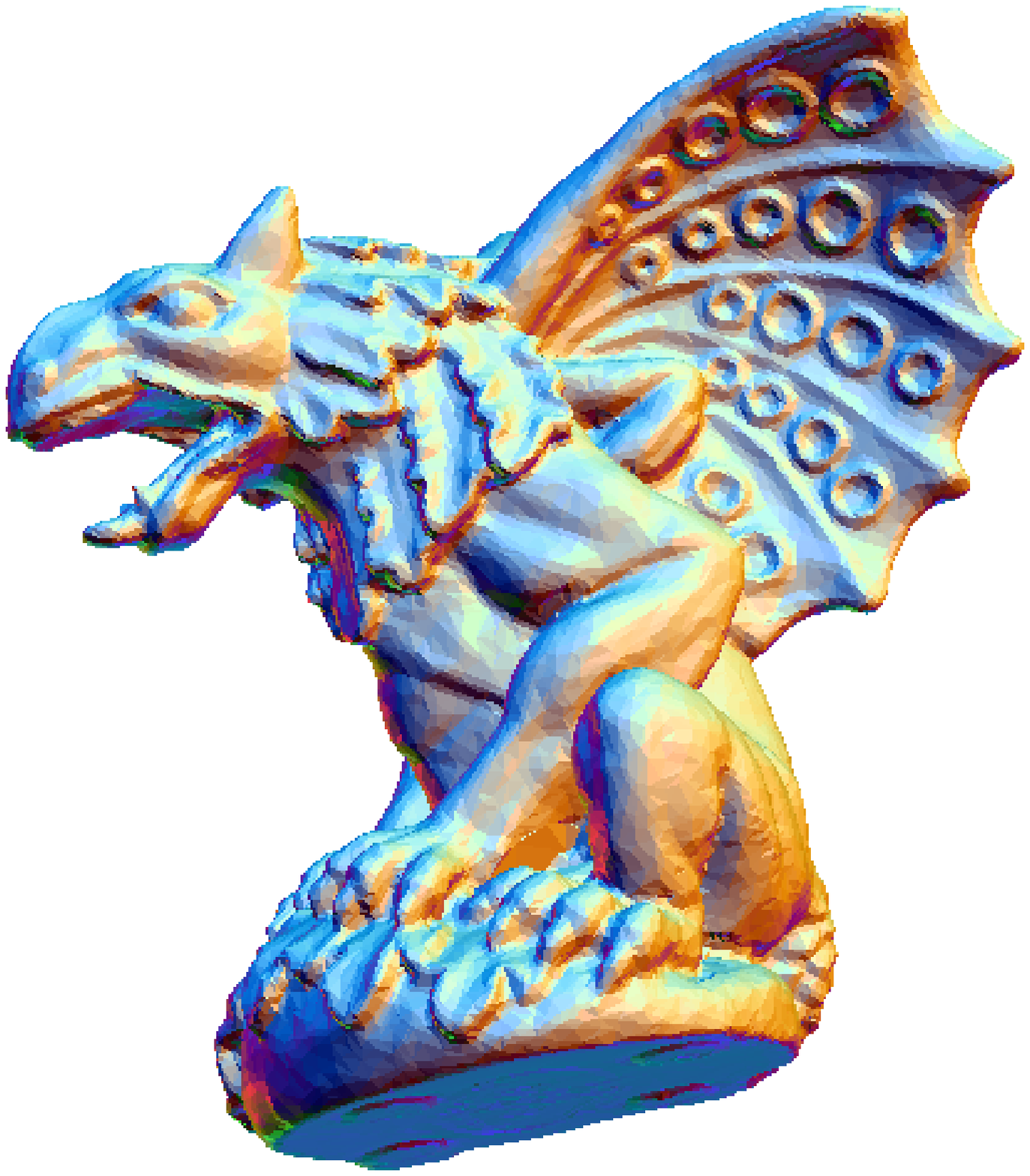}
\includegraphics[width=0.32\textwidth]{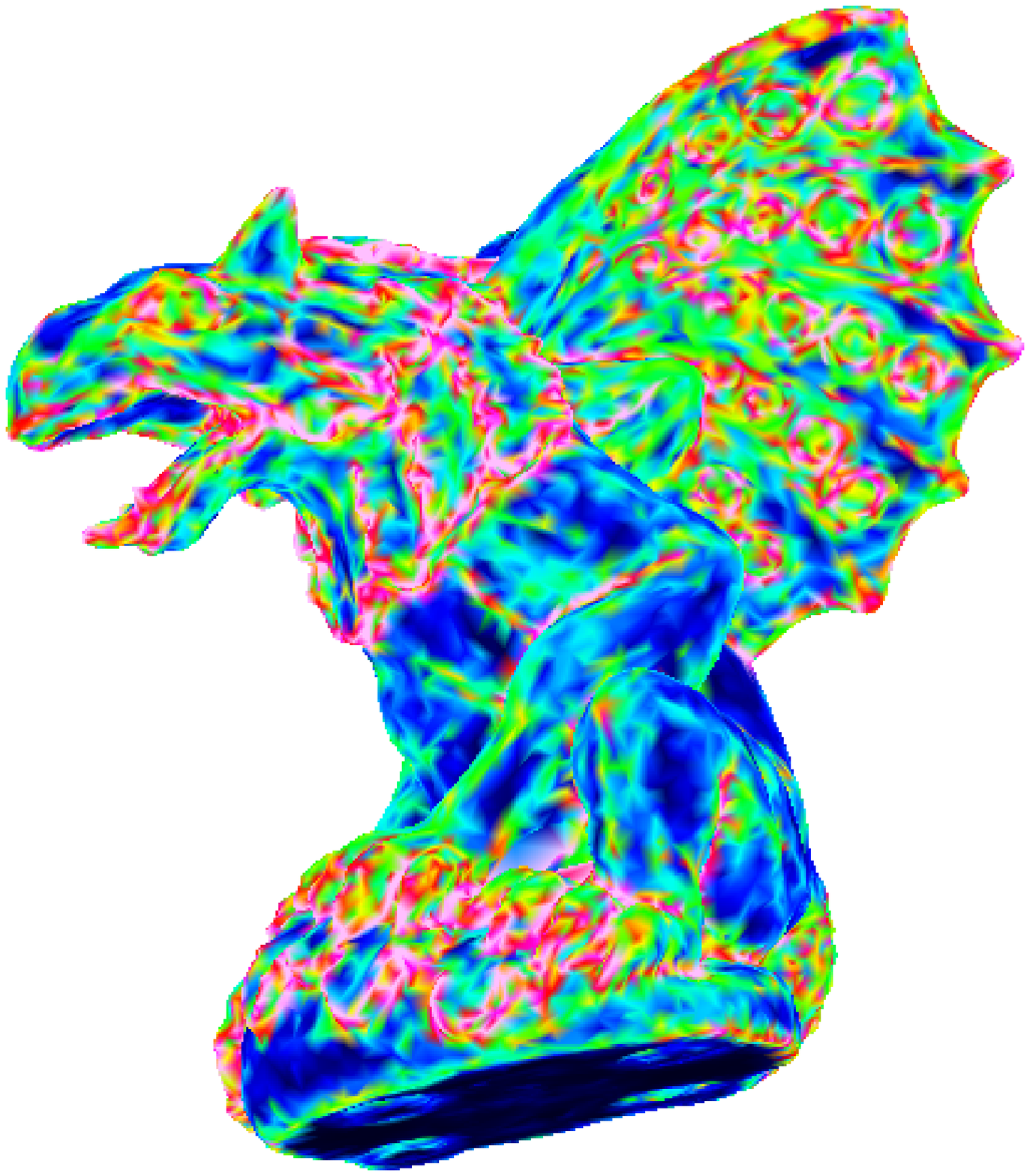}
\caption{Color encodings (please see the online version
of this paper for this figure in color)
for the mean curvature (left) and curvedness (right)
detected on Gargoyle model (center).
}
\label{HR-map}
\end{figure}

Figure~\ref{HR-map} presents color encodings for the mean curvature $H$
and curvedness $R$ detected on a geometrical model with many small-scale details (see the electronic version of this paper
for a color version of this figure).

A visual comparison of (\ref{HR-discrete}) and (\ref{HR-approx})
with the curvature estimation methods developed in \cite{Taubin_iccv95}
and \cite{pottmann-2009-iir} indicates that our discretization is less robust
to noise but better reflects small-scale surface details. One can improve the approximation
by using more sophisticated approximations for
the tangent directions, directional curvatures, and  directional derivatives
of the normal instead of common-sense estimates (\ref{HR-approx}).
Here our goal was to reveal potential applications of our approach
to curvature-based shape analysis.

\section{Discrete spherical means of derivatives}
\label{sec-t}

\paragraph*{\bf General construction and Veronese maps}

In this section we describe a generalization of the discrete mean value formulas of the preceding section from the circle case to arbitrary higher dimensions.
In particular, we prove the following
%
\begin{theorem} \label{main_problem}
%
Given a degree $k\not=0$ and a collection of points $P_i, i=1,\dots,N$, in general
position on the unit sphere in $S^n\subset \R^{n+1}$, with
$N > {n+k\choose n}-{n+k-2\choose n}$, there exist weights $w_i$
(not all equal to zero), which can be found constructively, 
such that for every harmonic homogeneous polynomial $F$
in $\R^{n+1}$ of  degree $k$ one has $\sum_i w_i \cdot F(P_i)=0$.
If $N > {n+k\choose n}$,  the same statement holds  for any
homogeneous polynomial $F$ in $\R^{n+1}$ of degree $k$.
\end{theorem}

For $k=1$, we have $N>n$, and this theorem can be viewed as a direct corollary of the Minkowski  theorem.
Note that for a possible application to computations of the mean curvature of a manifold, the function $F$ can be the sum of a constant ($k=0$) and a homogeneous quadratic term ($k=2$), and one needs to find the weights for obtaining the mean value $\bar F$ of such a function from its values at the given set of points. (Of course, all nonzero harmonics should give us the zero mean.) The cases of computations for the discrete Laplacian and square gradient operators are similar.

For instance,
{we use Theorem~\ref{main_problem}
to construct a discrete Laplacian as follows.}
{
\begin{definition}
{\rm
Let $O$ be a vertex in an arbitrary irregular grid with neighbouring vertices $\bar P_i$.
Define the discrete Laplacian of a function $f$ at a point $O$ to be
\begin{equation}\label{multi-Laplace}
\dfrac12\Delta f (O):= \sum w_i\dfrac{\partial^2 f}{\partial\e_{P_i}^2}\,,
\end{equation}
where points $P_i$ are normalized from $\bar P_i$ to lie on the unit sphere centered at $O$, the weights $w_i$ are given by {the above theorem for $k=2$}, and ${\partial^2f}/{\partial\e_{P_i}^2}$ are discrete approximations
of the second derivatives  of $f$  at $O$ in the direction $OP_i$.
}
\end{definition}
}

{
Definition~(\ref{multi-Laplace})
is inspired by formula (\ref{dmv4lap})
and, as we will see later, generalizes the standard discrete Laplacians
defined on regular grids.
}

\begin{remark}
{\rm
{The result of Theorem~\ref{main_problem}}
is very close to constructions of so called ``spherical designs." The latter
is a finite set of points on $S^n$ such that the average value of any polynomial $F$ of degree less than or equal to $k$ on this set equals the average value of the polynomial $F$ on the sphere. The concept of a spherical design was introduced in \cite{DGS}. The existence and structure of spherical designs in a circle was studied in \cite{Hong-1982}. In \cite{SZ}  the existence of designs of all sufficiently large sizes was proved: there is a number $N(n,k)$ depending on $n$ and $k$,  such that a design exists for any $N>N(n,k)$.
Many examples and recent results can be found in the survey \cite{Bannai}.

{The statement of Theorem~\ref{main_problem}}
is somewhat different: we do not require the weights to be equal.
In many applications, however, one uses the regular latices, and the weights turn out to be equal, as we discuss below.
}
\end{remark}

\paragraph*{\bf Proof--Construction}
%
First, we prove Theorem~\ref{main_problem} for a linear functional $L$
(i.e. for $k=1$). In this case the solution is given by the Minkowski
theorem~\ref{Mink-exist}. Namely, consider the tangent hyperplanes to
the sphere $S^n$ at points $P_i$; they form a polyhedron. Here we use
the fact that $P_i, i=1,\dots,N$, are in general position, so that no hyperplanes
are parallel and no more than $n+1$ intersect at a point (note that the polyhedron
is not necessarily circumscribed about the sphere, but the hyperplanes containing
its facets are tangent to the sphere). Let $w_i$ be the signed codimension-one
volumes of the faces of this polyhedron. Then $\sum_i w_i \cdot P_i =0$ according
to the Minkowski theorem (where we understand $P_i$ as vectors in $\R^{n+1}$ and
where the weights may be negative), and hence $\sum_i w_i \cdot L(P_i) = 0$ for
every linear function.


Next, consider  harmonic polynomials of degree $k$ in $\R^{n+1}$. Their restrictions
on the sphere form the space ${\cal H}_k$ of  spherical harmonics of degree $k$.
Consider the harmonic Veronese map $V: \R^{n+1} \to \R^q$, where
$V(x_1,...,x_{n+1})=(x_1^k-x_2^k,...)$, that is, the map defined by taking
a basis of harmonic polynomials in $x_1,...,x_{n+1}$ of degree $k$.
The dimension  is given by the formula:
$q=\dim {\cal H}_k = {n+k\choose n}-{n+k-2\choose n}$.
Then $F=L \circ V$ where $L$ is a linear function on $\R^q$. Let $w_j$
be the constructed above weights for the collection of points $V(P_j)$ in $\R^q$.
This gives us the required formula $\sum_i w_i \cdot F(P_i) =0$. Note that
the general position condition on $P_i$ now means that the points
$Q_i=V(P_i)$ in $R^q$, after normalization, satisfy the above condition that
the tangent planes to the sphere $S^{q-1}$ form a polyhedron.


Finally consider all homogeneous polynomials in $\R^{n+1}$ of  degree $k$.
Their restrictions on the sphere form the space
${\cal H}_k\oplus {\cal H}_{k-2}\oplus {\cal H}_{k-4}\oplus\ldots$.
The dimension of this space is ${n+k\choose n}$, and we repeat
the same construction involving the Veronese map for each homogeneous
component as before.

\begin{example}\label{example2d}
{\rm
Let us demonstrate this approach in the case of a circle $S^1$ and quadratic functions.
The harmonic Veronese map for quadratic polynomials sends the coordinate functions $\{x_1, \dots, x_{n+1}\}$ to the space of quadratic harmonic polynomials (whose restrictions on the sphere $S^n$ give us spherical harmonics). For the circle $S^1\subset \R^2$ this harmonic Veronese map is
$V: \R^2\to \R^2,$ where $(x,y)  \mapsto (x^2-y^2, 2xy)$, or $(\cos \phi, \sin \phi) \mapsto (\cos 2\phi, \sin 2\phi)$, see a detailed analysis in the preceding section.

Note that if we map $\R^2_{x,y}$ to the full 3-dimensional space $\R^3$ of quadratic polynomials
$(x^2, y^2, 2xy)$, then the  image of the circle $S^1\subset \R^2$
will lie in the plane $\R^2\subset \R^3$(due to the relation $x^2+y^2=1$),
so the point images $V(P_i)$ will not be in general position to yield a polyhedron circumscribed around a 2-sphere in $\R^3$.
}
\end{example}

\begin{example}
{\rm
For the case of 2-sphere and quadratic functions  the space of harmonic polynomials is
5-dimensional, and our harmonic Veronese map sends $S^2$ to this $\R^5$. The point images lie in general position, and sufficiently many points (at least 6) would yield a polyhedron, which allows one to find the corresponding weights $w_i$.
}
\end{example}

\begin{remark}
{\rm Since the points $Q_j=V(P_j)$ should lie on the sphere $S^{q-1}$, rather than be generic vectors in $\R^q$, one can project them to the sphere by rescaling. Namely, let $Q_j$ be a (sufficiently generic) collection of non-zero vectors in $\R^q$. Assign the following weights $w_j$ to these
vectors. Consider the collection of points $Q_j/|Q_j|$ on the unit sphere $S^{q-1}$, draw the tangent hyperplanes to the sphere, and let $u_j$ be the signed $(q-1)$-dimensional volumes of the faces of the resulting polyhedron.
}
\end{remark}

\begin{proposition}\label{Prop-key}
For every linear function $L$ on $\R^q$, one has:
$\sum_j w_j\cdot L(Q_j)=0$, where the weights are $w_j:=u_j/|Q_j|$.
\end{proposition}

\noindent
{\sc Proof.} The Minkowski formula says that $\sum u_j \cdot Q_j/|Q_j| =0$, hence $0=\sum u_j \cdot L(Q_j)/|Q_j| =\sum_j w_j\cdot  L(Q_j)$. QED

\bigskip

\paragraph*{\bf Weights for regular lattices}

Now, from the general consideration of points in general position we move to regular ones,
and consider the set of all ``$m$-mid-points,"
i.e., the set $M_m\subset S^{n}$ of points  whose coordinates are all possible combinations
of $m$  nonzero coordinates equal to $\pm1$ among $n$ spots, starting with $(\pm1, ...,  \pm1, 0,...,0)$ to $(0,..., 0,\pm1, ...,  \pm1)$, which are scaled to belong to the unit sphere $S^n\subset \R^{n+1}$. Here $m$ is any fixed integer between 1 and $n$. We also confine to the case of quadratic harmonics, $k=2$ in the general problem above.

\begin{example}
{\rm
In dimension 3, for $m=1$ we get $M_1$ consisting of 6 midpoints of the cube faces, which are also the vertices of an octahedron:
$(\pm1, 0,0), (0,\pm1,0)$ and $(0,0,\pm1)$. The set $M_2$ consists of 12 scaled mid-points of the cube edges:
$(\pm1, \pm1, 0)/\sqrt2$, $(\pm1, 0,\pm1)/\sqrt2,$ and $(0,\pm1, \pm1)/\sqrt2$. The set $M_3$
contains exactly 8 vertices of the cube: $(\pm1, \pm1, \pm1)/\sqrt3$.
}
\end{example}

\begin{theorem}\label{theor}
Let $F(x_1,\dots,x_{n+1})$ be a homogeneous quadratic function, $M_m\subset S^{n}$ be the set of
$m$-mid-points for any fixed $m$, $1\le m\le n$.  Then
$$
\frac{1}{vol(S^{n})} \int_{S^{n}} F~ dS=\frac{1}{\#M_m}\sum_{x\in M_m} F(x)\,,
$$
where $vol(S^{n})$ is the volume of $n$-dimensional sphere and $\#M_m$ is the number of points in $M_m$.
\end{theorem}

One can see that in all these ``regular cases" one counts the points of $M_m$ with equal weights.
This theorem generalizes the corresponding 3D consideration
(see e.g. \cite{Patra-Karttunen_nmpde05}), which gives three basis stencils used in the Laplacian computations.

\noindent
{\sc Proof.} Write $F=c+G$, where $c\in \R$ and $G$ a harmonic  quadratic polynomial. For a constant $c$ the claim is obvious. For $G$ one has $\int_{S^{n}} G dS=0$, so we need to show that $\sum_{x\in M_m}G(x)=0$ for every harmonic quadratic polynomial $G$.

Consider the full Veronese map $\overline V: \R^n \to \R^{q+1}$, given by
$$
(x_1, ..., x_{n+1})\mapsto (x_1^2, ... , x_{n+1}^2; x_1 x_2, ..., x_{n}x_{n+1})\,,
$$
with $q+1=(n+1)(n+2)/2$.
Let $\cal H$ be the space of  linear functions on $\R^{q+1}$ that vanish on the vector
$\xi=(1,...,1;0,...,0)$ with first $n+1$ nonzero components (we separate them by semicolon).
Then, for each quadratic harmonic $G$, one has $G=H\circ \overline V$ for some $H\in \cal H$. Indeed, harmonic polynomials form a hyperplane in the space of linear functions on $\R^{q+1}$ being  given by one linear relation (coming from the sphere equation $x_1^2+...+x_{n+1}^2=1$), and the polynomials $x_1^2-x_2^2, ..., x_n^2-x_{n+1}^2, x_1 x_2, ..., x_nx_{n+1}$ are harmonic and equal to 0 on this vector $\xi$. Thus we need to check that $\sum_{Q\in \overline V(M_m)}H(Q)=0$ for all $H\in \cal H$ and any $m=1,...,n+1$.

For instance, for the set $M_1$, which consists of $(0,...,0,\pm1,0,...,0)$-type points, the image  $\overline V(M_1)$ consists of vectors $Q_i$:
$ (0,...,0,1,0,...0; 0,  .. 0)$ with 1's in one of the first $n$ coordinates.
In the case of $M_2$ with $2n(n-1)$ ``diagonal" points
$$
(0,...,0,\pm1, 0,..., 0, \pm1, 0,...,0)/\sqrt 2,
$$
the image $\overline V(M_2)$ consists of vectors $Q_i$:
$$
(0,...,0,1,0,...0, 1,0,..., 0\,; 0, ...,0, \pm1, 0, ..., 0)/2
$$
with all possible pairs of 1's among first $n+1$ coordinates and the $\pm1$ at the respective, $ij$-place. Similarly, for any $m$ one can see that in each of these cases we have  $\sum_1^{\#M_m} Q_i=const\cdot \xi$. Hence    $\sum_1^{\#M_m} H(Q_i)=0$ for all $H\in \cal H$, as claimed.
QED.

Thus, while Theorem~\ref{main_problem} allows one to find the weights $w_i$ for a generic set of vectors $P_i$,  Theorem~\ref{theor} gives the explicit weights in the case when $F$ is a homogeneous quadratic polynomial (i.e., $k=2$) and vectors $P_i$ correspond to the $m$-mid-points of the $n$-dimensional cube, $1\leq m\leq n$.

\begin{remark}
{\rm
 For computations we are often interested mostly in the case $S^2\subset \R^3$.
Here is the explicit statement for $M=M_2$ adjusted to this case:
let $F(x,y,z)$ be a homogeneous quadratic function, $M\subset S^2$ be the set of 12 points:
$(\pm1, \pm1, 0)/\sqrt 2$, $ (\pm1, 0, \pm1)/\sqrt 2$, and $ (0, \pm1, \pm1)/\sqrt 2  $.
Then
$$
\frac{1}{4\pi} \int_{S^2}F~dS=\frac{1}{12}\sum_{P\in M} F(P).
$$

Thus in  all of the above cases we need to use the same weights on these points.
This equality of weights for a given set $M_m$ can also be obtained from the symmetry consideration. One can also take combination of the above lattices $M_m$ with
with different weights for different $m$
(provided that the total sum of the weights is equal to 1),
e.g., by considering arbitrary combinations of the 3 lattices in $\R^3$: mid-edges,
vertices of the cube, and at the middle of faces. By imposing additional requirements, one can optimize these coefficients: for instance, to achieve an isotropic approximation, or with the next approximation smallest in an appropriate sense.
}
\end{remark}

\begin{remark}
{\rm
Points of $\overline V(M_1)$, evidently, form the coordinate $(n+1)$-simplex in the space spanned by first $n+1$ coordinates, which is a slim subspace in $\R^{q+1}$. In the case of $M_2$, the
points of $\overline V(M_2)$ form the standard $q$-simplex in a hyperplane in $\R^{q+1}$.
Indeed, in the latter case, we see that
for all $Q_i\in \overline V(M_2)$ one has $(Q_i, Q_j)=1/2$ for $i\not=j$ and $(Q_i, Q_i)=3/4$, i.e., the vectors $Q_i$ have the same length and same angles between them.
}
\end{remark}

\medskip

\begin{remark}
{\rm
Above we considered two discretization procedures, both using the Veronese map: a generic set of points $\{P_i\}$ and the mid-points sets. While for a generic set in $\R^{q+1}$ the coefficients
are given by the Minkowski formula for the corresponding polyhedron in the image, it is not the case in the mid-point cases. Indeed, due to a  very regular structure of the image points $Q_j=\overline V(P_j)$, the corresponding tangent planes go to infinity, and the corresponding volumes of the faces are infinite.

In a sense, if the image points $Q_j$ are linearly dependent in $\R^{q+1}(\mod \xi)$,
and the corresponding polyhedron has faces going to infinity (like in the standard basis case), we have to quotient out these ``directions to infinity," and consider the volume of faces in a polyhedron of smaller dimension.
}
\end{remark}

\begin{example}
{\rm
For $n=2$ and for the standard basis we did not get a polygon on the plane $\R^2= \R^3 (\mod \xi)$. Instead, it reduces to a pair of points in $\R^1$.
It is not enough to start with two generic points in this case, but if we take three generic points on the circle, we obtain an interpolation formula.
This is a manifestation of the general phenomenon described in the following proposition.
}
\end{example}

\begin{proposition}
There is no discrete mean value formula for quadratic harmonics and a given generic collection of $n+1$ points on $S^{n}$.
\end{proposition}

\noindent
{\sc Proof.}
Let $P_1,...,P_{n+1} \in S^{n}$, and let $\overline V:\R^{n+1} \to \R^{q+1}$ with $q+1=(n+1)(n+2)/2$
be the full Veronese map:
$$(x_1,...,x_n) \mapsto (x_1^2,...,x_{n+1}^2;x_1x_2,...,x_ix_j,...).$$
Let $\overline V(P_k)=Q_k$. Let $\xi=(1,...,1;0,...,0)\in \R^{q+1}$.

We want to find constants $w_k$  such that $\sum w_k \cdot Q_k=\xi$. This is equivalent to the following equality:
$$
Diag(w_k) =  (A^{*} A)^{-1},
$$
where $A$ is the matrix made of the vectors $P_1,...,P_{n+1}$, as the following linear algebra shows.
Indeed, the equality $\sum w_k Q_k=\xi$ is equivalent to
$$
\sum_i w_i x_{k,i} x_{l,i}=\delta_{k,l},
$$
which is equivalent to $A\cdot Diag(w)\cdot A^{*} =E$, the identity matrix. This is equivalent to $(A^{*} A)^{-1}=Diag(w)$, as claimed.

Of course, if $\{P_k\}$ is an orthonormal frame then $A$ is orthogonal, and all $w_k=1$, as in the theorem above.
The matrix $AA^{*}$ is diagonal if and only if the vectors $P_i$ are pairwise orthogonal.
This is the case for the set $M_1$ discussed above, which consists of the basis points, but, of course, not true for generic collection of $n+1$ points on $S^{n}$. This  implies that
one cannot have an interpolation formula with $n+1$ generic points. QED

Note that on the other hand, if one has sufficiently many points, so that their Veronese images are linearly dependent in $\R^{n(n+1)/2}\, (\mod \xi)$, then any such linear relation provides the weights that yield an interpolation formula.


\paragraph*{\bf Higher-dimensional Laplacian approximation formulas}

The results of the present section allow one to easily construct various discrete Laplacians and squared gradient  similar to representations (\ref{dmv4lap}) and (\ref{dmv4grad}), respectively.
In the beginning of the section we discussed the construction of the discrete Laplacian for an arbitrary grid, see formula (\ref{multi-Laplace}):
$$
\dfrac12\Delta f (O)= \sum w_i\dfrac{\partial^2 f}{\partial\e_{P_i}^2}\,.
$$
To obtain the weights $w_i$ for a point $O$ with neigboring vertices $\bar P_i$ we first normalize the latter to get $P_i$ lying on the unit sphere centered at $O$. Then we
consider the Veronese map and obtain points $Q_i$. Finally we obtain the weights $w_i=u_i/|Q_i|$ by finding the signed $(q-1)$-dimensional volumes $u_i$ of the faces of the polyhedron in the $Q$-space, see Proposition \ref{Prop-key}. An analog of the formula (\ref{dmv4grad}) for the square gradient instead of the Laplacian  is completely analogous.

\medskip

For a regular grid, Theorem \ref{theor} delivers similar formulas for a discrete Laplace operator where all weights $w_i$ are equal:
\begin{equation}\label{Laplace-m}
\Delta f (O)= C_m\sum \dfrac{\partial^2 f}{\partial\e_{P_i}^2}\,,
\end{equation}
where points $P_i$ belong to the grid of $m$-midpoints, $P_i\in M_m$, and the constant $C_m$ normalizes the sum, $C_m=\#M_1/\#M_m$, by comparing the number
of neighbours in $M_m$  with that in the coordinate grid $M_1$, the most straightforward definition of the discrete Laplacian.

Furthermore, one can consider  linear combination of the Laplace approximations (\ref{Laplace-m}) with different coefficients $\alpha_m$ satisfying $\sum_{m=1}^n \alpha_m=1$.
Any such linear combination delivers a different  discrete Laplacian and one can run a separate optimization problem on these coefficients $\alpha_m$ to find
a discrete Laplacian with better rotation-invariant properties. We discuss the work in this direction for the 2D Laplacian in the next section.


\section{\bf Discrete circular means of derivatives and isotropic
finite differences}\label{sec-l}
In this section we show how to use discrete circular means of Section \ref{sec-a} 
for approximation of various differential operators in 2D by finite differences
with good rotation invariance properties.

The simplest way to construct
such approximations consists of substituting in (\ref{dmv4lap}), (\ref{dmv4grad}), (\ref{dmv4dirder}), and (\ref{qLap}) three-point finite
differences instead of the first-{} and second-order derivatives.
For example, constructing a discrete Laplacian corresponding to
(\ref{dmv4lap}) can be done as follows. Consider a point $\p$ inside
a convex bounded domain $\Omega$. Let the straight line determined by $\p$
and direction $\e_\varphi$ intersect the boundary $\partial\Omega$
in two points $\q_1$ and $\q_2$. Denote by $\rho_1$ and $\rho_2$
the distances from $\p$ to $\q_1$ and $\q_2$, respectively.
The second-order $\e_\varphi$-directional derivative of $f(\p)$
can be approximated by
\begin{equation}\label{fd2}
\dfrac{\partial^2 f}{\partial\e_\varphi^2}\approx
\dfrac{2}{\rho_1\rho_2}
\left(\left.\left[\dfrac{f(\q_1)}{\rho_1}
+\frac{f(\q_2)}{\rho_2}\right]\right/
\left[\frac{1}{\rho_1}+\frac{1}{\rho_2}\right]-f(\p)\right),
\end{equation}
where $\rho_1$ and $\rho_2$ are assumed to be small.
Now a discrete approximation of the Laplacian is obtained
by summing up (\ref{fd2}) over the set of directions
$\e_{\varphi_k}$, $k=1,2,\dots,n$.

Discrete Laplacian (\ref{dmv4lap}), (\ref{fd2}) can be considered
as a special case of a general construction studied in
\cite{Wardetzky_sgp07} (see also references therein).

It is interesting to note that integration of (\ref{fd2})
with respect to $\varphi$
\begin{equation}\label{wGW}
\Delta f(\p)
\approx\left.\int\limits_0^{2\pi}\left[\left.\left(\frac{f(\q_1)}{\rho_1}
+\frac{f(\q_2)}{\rho_2}\right)\right/
\left(\frac{1}{\rho_1}+\frac{1}{\rho_2}\right)\right]
\,\frac{d\varphi}{\rho_1\rho_2}
\right/
\int\limits_0^{2\pi}\frac{d\varphi}{\rho_1\rho_2}-f(\p)
\end{equation}
leads to an extension of a pseudo-harmonic interpolation
scheme proposed in \cite{Gordon-Wixom_sinum74}. It can be shown
that, given $f(\q)$ defined on $\partial\Omega$, (\ref{wGW}) and its generalizations deliver the harmonic interpolation if $\Omega$ is
a circle (a ball in the multidimensional case) \cite{Belyaev_sgp06}.

Unfortunately combining circular means with directional finite
differences does not automatically guarantee rotation invariance
properties of the approximating differential operators. The reason
is that the finite differences introduce directional errors. However,
as seen below, for a regular grid, such directional errors
can be uniformly distributed over the directions by choosing
appropriate weights in discrete circular means (\ref{dmv4lap}),
(\ref{dmv4grad}), (\ref{dmv4dirder}), and (\ref{qLap}).

\paragraph*{\bf Gradient and Laplacian stencils for square grids}
Now let us focus on analyzing rotation invariance properties of
regular grid stencils corresponding to (\ref{dmv4lap}) and
(\ref{dmv4grad}).

Consider a two-dimensional square grid with step-size $h\ll1$.
Each grid vertex has eight nearest neighbors: two horizontal, two vertical,
and four diagonal.
Thus, according to our geometric interpretation of
(\ref{dmv4lap}), (\ref{dmv4grad}), and (\ref{dmv4dirder}), we are given
four directions and corresponding unit normals
$$
\left\{e^{2i\varphi}\right\},\quad \varphi=0,\pi/4,\pi/2,3\pi/4
$$
which define a rectangle aligned with the coordinate axes.

Consider such a rectangle with the edge lengths
$2$ and $w$, as shown in the right image of Figure~\ref{Mink},
and replace the first-{} and second-order directional derivatives
in (\ref{dmv4lap}), (\ref{dmv4grad}), and (\ref{dmv4dirder})
by corresponding central differences.
We arrive at the following stencils for the $x$-derivative and Laplacian
\begin{equation}\label{wgrad}
\frac{\partial}{\partial x}\approx D_x\equiv
\frac{1}{2h(w+2)}\left[\begin{array}{rcl}
-1 & 0 & 1
\\
-w & 0 & w
\\
-1 & 0 & 1
\end{array}\right]
\end{equation}
and
\begin{equation}\label{wLap}
\Delta\approx\wLap\equiv
\frac{1}{h^2(w+2)}\left[\begin{array}{rcl}
1 & w & 1
\\
w & -4(w+1) & w
\\
1 & w & 1
\end{array}\right]\,,
\end{equation}
respectively.
The stencil $D_y$ for $\partial/\partial y$ is obtained from $D_x$
by $\pi/2$-rotation. Here and everywhere below the $3\times3$ matrices are understood
as finite difference operators acting on functions defined on a square grid with
step-size $h$.

Formulas (\ref{wgrad}) and (\ref{wLap}) provide us with consistent parameterizations of the $3\times3$ stencils for the gradient and Laplacian.

\begin{remark}
{\rm The standard central difference and five-point Laplacian $\L_{+}$ is obtained for  $w=\infty$.
Setting $w=1$ leads to the so-called Prewitt masks for the 1st-order derivatives and
Laplacian \cite{Pratt_dip3e01,Gonzalez-Woods_dip3e08}. The case $w=2$ is commonly
used in image processing applications and yields the standard Sobel mask for the
derivative \cite{Pratt_dip3e01,Gonzalez-Woods_dip3e08} and a 9-point discrete Laplacian,
which possesses good isotropic properties \cite{Kamgar-Parsi_tip99}.
The case $w=4$ was analyzed in \cite{Bickley_qjmam48}, where it was shown that
\begin{eqnarray}
\label{grad4}
\Dx4 &=&
\dfrac{1}{8h}\left[
\begin{array}{rcl}
-1 & 0 & 1 \\
-4 & 0 & 4 \\
-1 & 0 & 1
\end{array}
\right]=\dfrac{\partial}{\partial x}+\dfrac{h^2}{12}\,
\Delta\dfrac{\partial}{\partial x}+O(h^4),
\\
\label{L4}
\Lap &=&
\dfrac{1}{6h^2}\left[
\begin{array}{rcl}
1 & 4 & 1 \\
4 & -20 & 4 \\
1 & 4 & 1
\end{array}
\right]=\Delta+\dfrac{h^2}{12}\,\Delta^2+O(h^4),
\end{eqnarray}
as $h\to0$.
Since the Laplacian is isotropic, the right hand sides of (\ref{grad4})
and (\ref{L4}) deliver asymptotically optimal asymptotically
rotation-equivariant $3\times3$ stencils for the $x$-derivative and Laplacian, respectively, as $h\ll1$.
In particular, this explains why (\ref{wLap}) with $w=4$ is often used
for a numerical solution of the Poisson equation
\cite[Chapter 3, \S1]{Kantorovich-Krylov_amha56},
\cite[Chapter 3, \S10]{Birkhoff-Lynch_nsep84}.
}
\end{remark}

Discrete nine-point Laplacian (\ref{wLap}) can be represented as a linear combination
of two basis five-point Laplacians
$$
\wLap=\alpha\L_{+}+\beta\L_{\times}
\mbox{ with }
\alpha=\dfrac{w}{w+2}\mbox{ and }\beta=\frac{2}{w+2},\mbox{ where}
$$
$$
\L_{+}=
\frac{1}{h^2}
\left[\begin{array}{rcl}
0 & 1 & 0
\\
1 & -4 & 1
\\
0 & 1 & 0
\end{array}\right]
\quad\mbox{and}\quad
\L_{\times}=
\frac{1}{2h^2}
\left[\begin{array}{rcl}
1 & 0 & 1
\\
0 & -4 & 0
\\
1 & 0 & 1
\end{array}\right].
$$
Setting $w=4$ yields $\alpha=2/3$ and $\beta=1/3$.

\begin{figure}[h]
\centering
\begin{tabular}{cccc}
\includegraphics[width=0.22\textwidth]{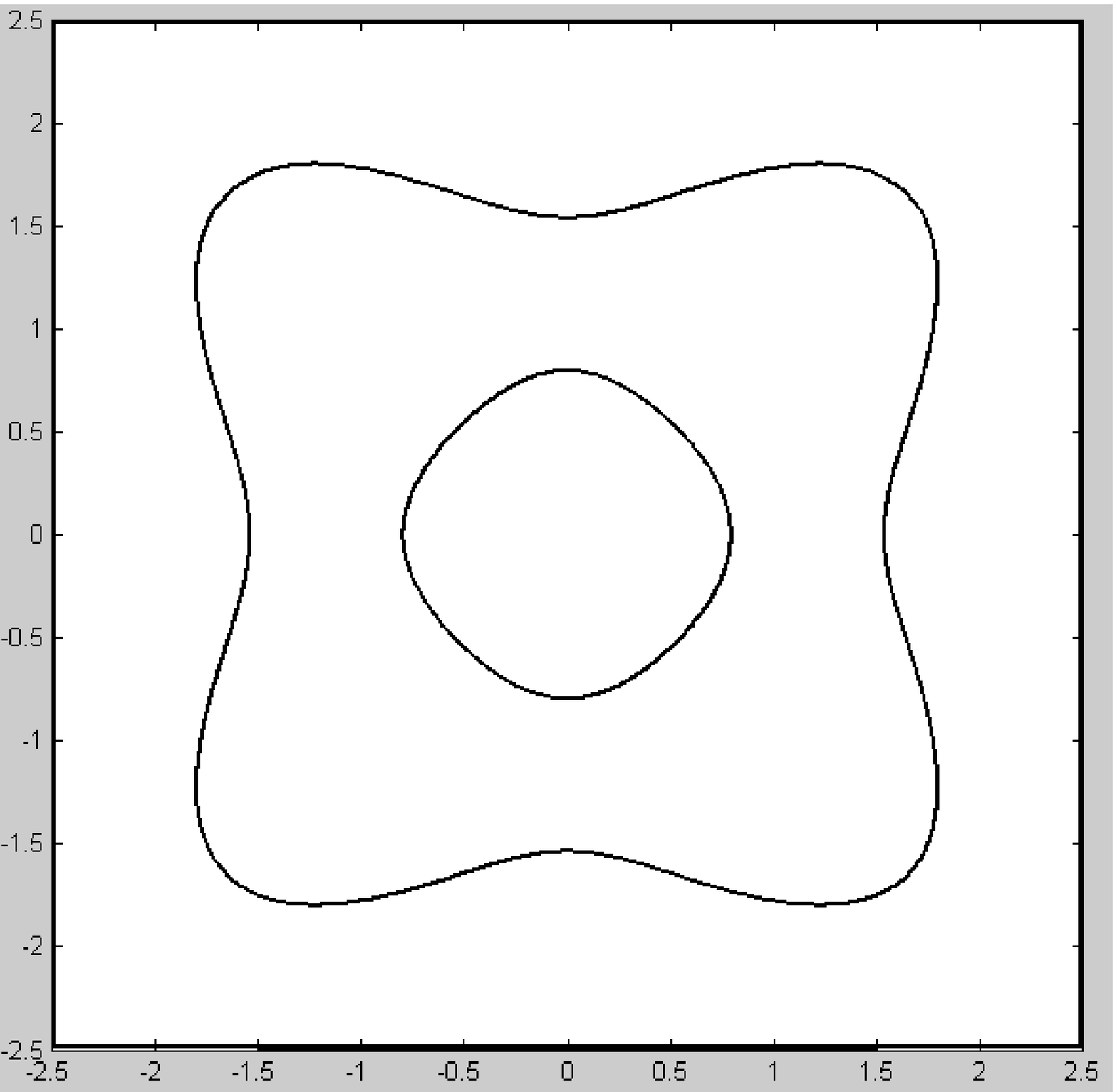}
&
\includegraphics[width=0.22\textwidth]{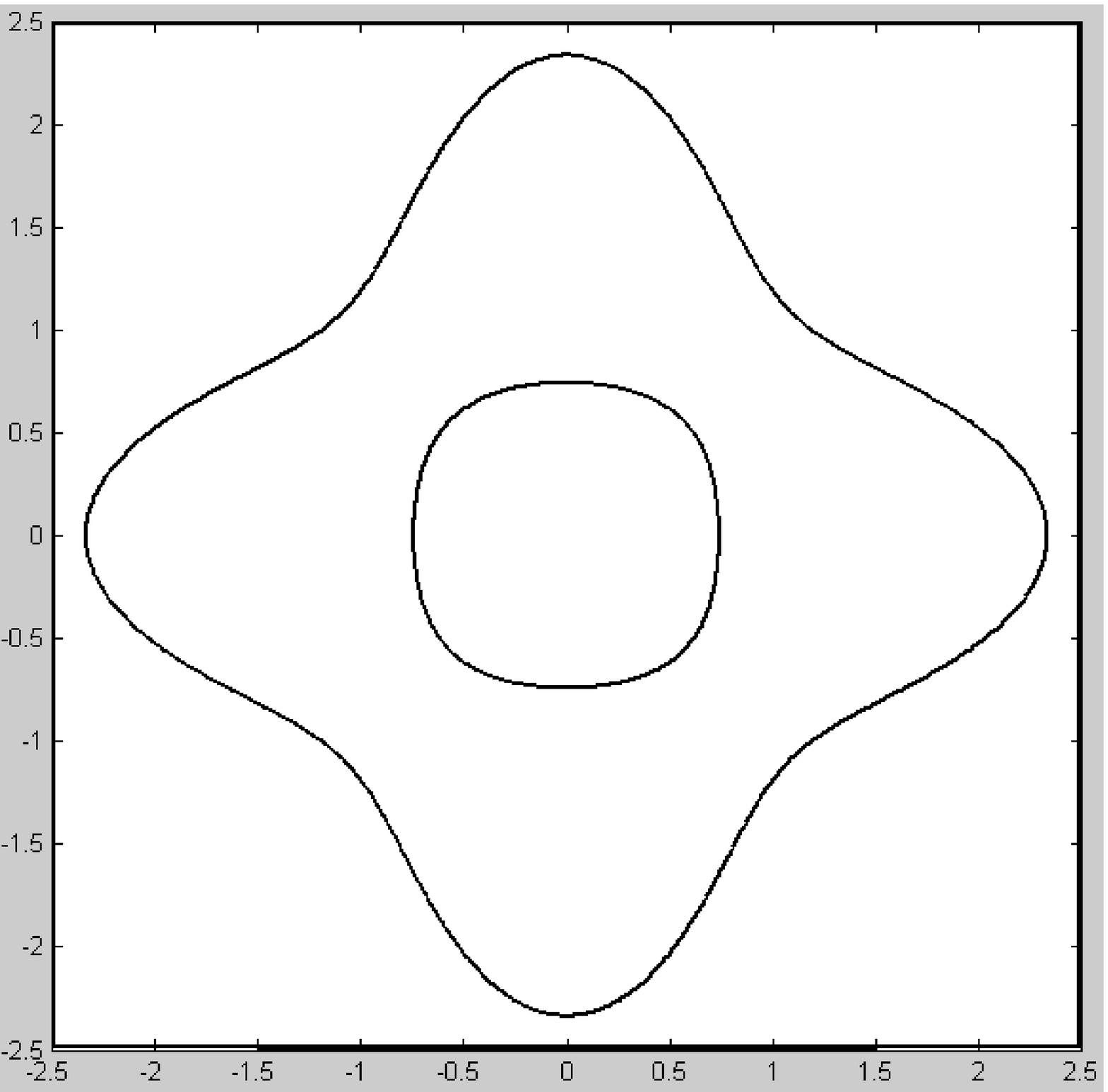}
&
\includegraphics[width=0.22\textwidth]{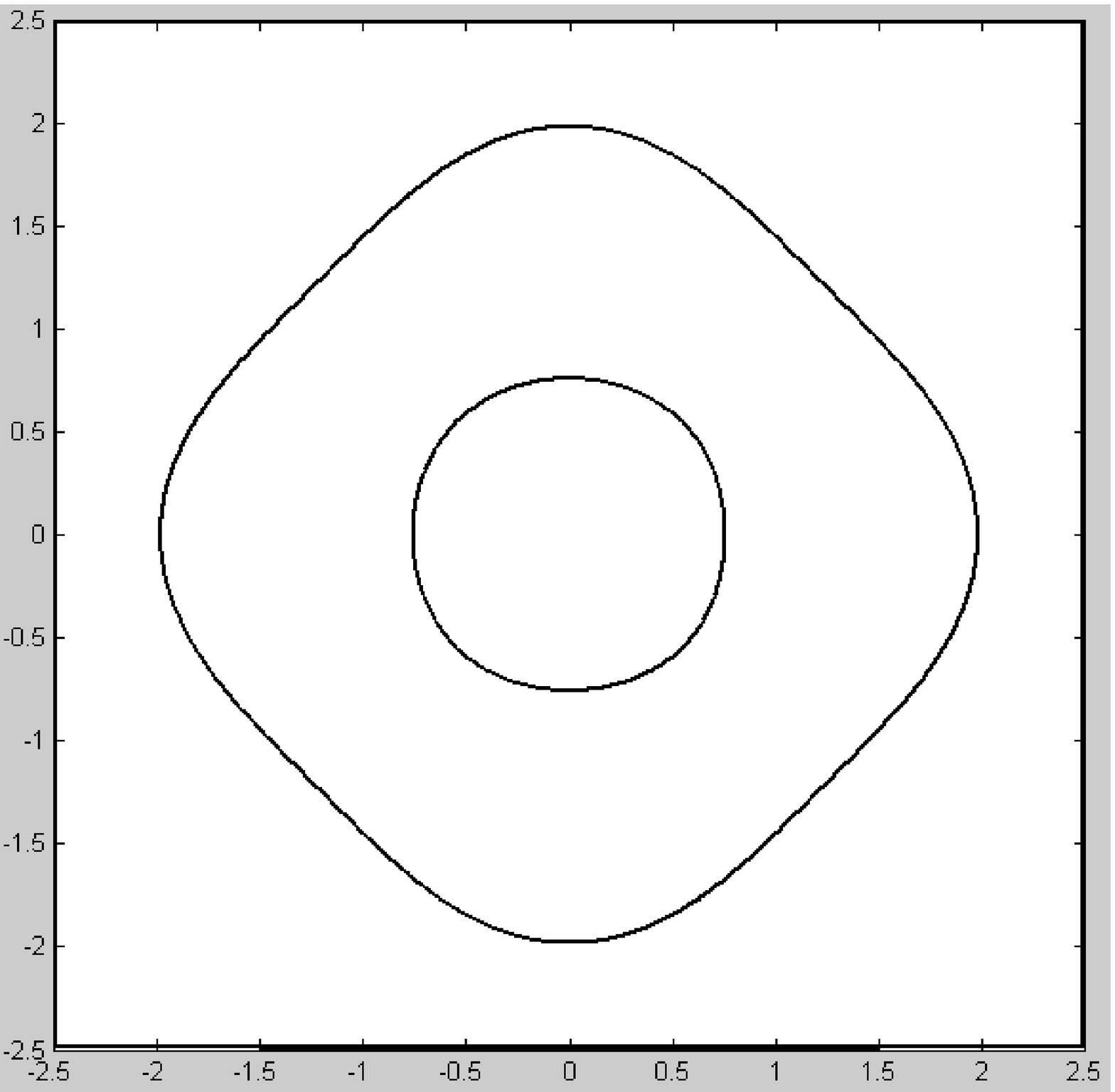}
&
\includegraphics[width=0.22\textwidth]{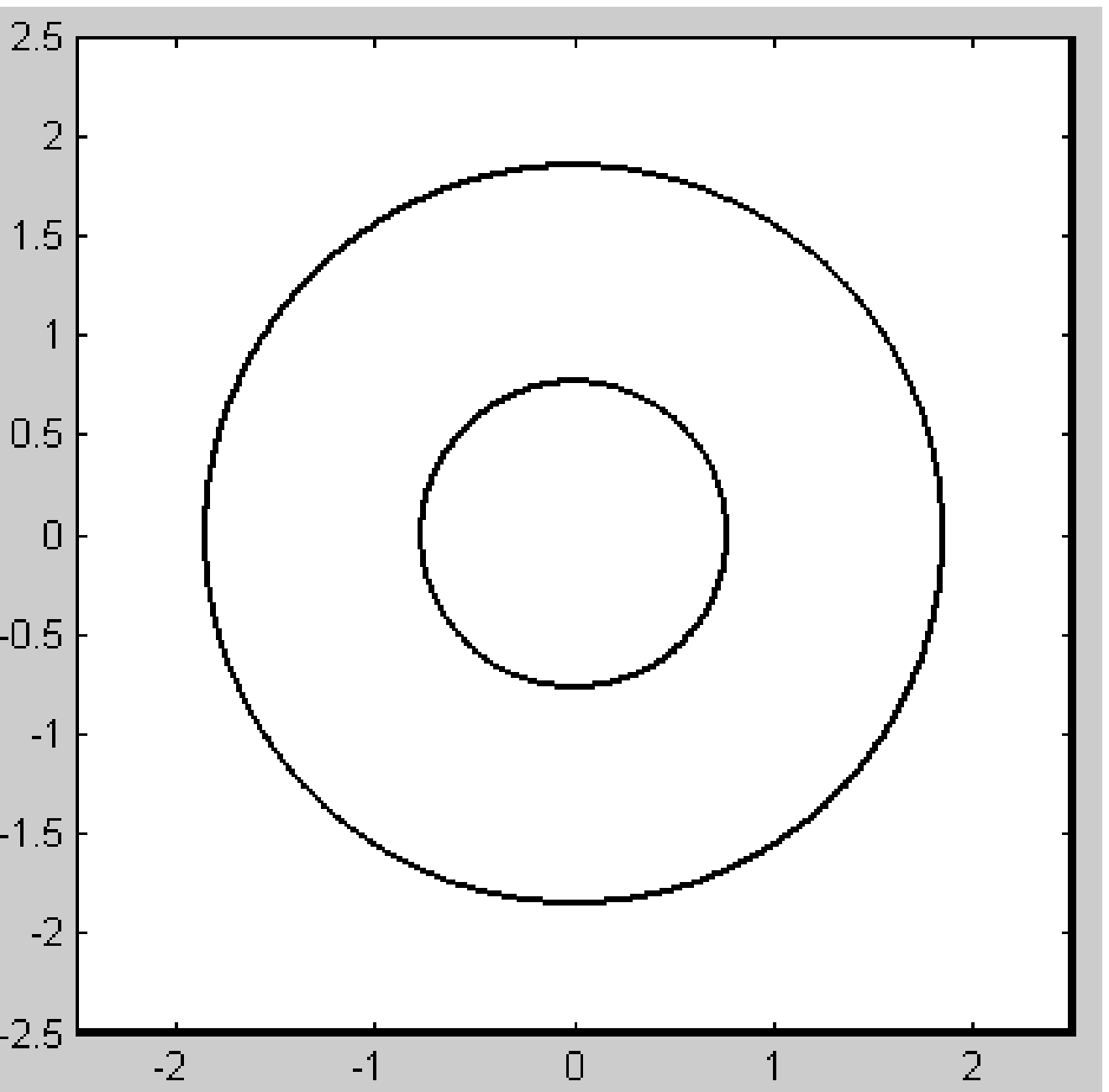}
\\
$w=\infty$
&
$w=0$
&
$w=2$
&
$w=4$
\\
{\small ($\alpha=1$, $\beta=0$)}
&
{\small ($\alpha=0$, $\beta=1$)}
&
{\small ($\alpha=1/2$, $\beta=1/2$)}
&
{\small ($\alpha=2/3$, $\beta=1/3$)}
\end{tabular}
\caption{Contour lines of $(\Delta-\wLap)[g]$ for various values of $w$.
The best performance is achieved for $w=4$.
}
\label{lap}
\end{figure}

In Figure~\ref{lap}, we use a Gaussian bump function
\begin{equation}\label{gaussian}
g(x,y)=\exp\left(x^2+y^2\right)
\end{equation}
to test isotropic properties of four discrete stencils
$$
\Lapi\equiv\L_{+},\quad\Lapz\equiv\L_{\times},
\quad\Lapt,\mbox{ and }\Lapf.
$$
Namely, the figure displays contour
lines of $(\Delta-\wLap)[g]$. As (\ref{L4}) suggests, $\wLap$ with $w=4$
demonstrates the best performance if the grid step-size $h$ is sufficiently
small.

Similar results are obtained for the $x$-derivative $3\times3$ stencils
(\ref{wgrad}) and corresponding $y$-derivative stencil. Analogously to the case
of the Laplacian, expansion (\ref{grad4}) implies that (\ref{wgrad}) with $w=4$
is rotationally optimal for sufficiently small $h$.

\begin{remark}
{\rm It is easy to increase the approximation accuracy of (\ref{grad4})
and (\ref{L4}). For example, one can rewrite (\ref{grad4}) as
$$
\frac{\partial}{\partial x}=
\left(1+\frac{1}{h^2}\Delta\right)^{\!-1}\!\!\!\Dx4+O(h^4)=
\left(1-\frac{1}{h^2}\Delta\right)\!\!\Dx4+O(h^4)
$$
and use stencil (\ref{L4}) for approximating $\Delta$.
}
\end{remark}


\paragraph*{\bf Quasi-Laplacian stencil with good rotation invariance properties}
Consider now a quasi-Laplacian operator
\begin{equation}\label{alap}
\L^a\equiv\nabla\cdot[a(x,y)\nabla]=a\Delta+\nabla a\cdot\nabla
\end{equation}
where $a(x,y)$ is given. In view of (\ref{L4}), (\ref{grad4}),
and (\ref{qLap}), the optimal $3\times3$ stencils for $\L^a$
are obtained from (\ref{wLap}) and (\ref{wgrad}) when $w=4$.

In practice, the diffusion coefficient $a(x,y)$ in (\ref{alap})
is defined at the staggered grid points which lie on halfway between the standard
grid points. Following
\cite[A.3.2]{Aubert-Kornprobst_mpip2e06} let us approximate
(\ref{alap}) using a linear combination of two basis stencils
$$
\alpha \L^a_{+}+\beta \L^a_{\times}
$$
with
$$
\begin{array}{l}
\L_{+}^a[f]\equiv\dfrac{1}{h^2}\left[
a_{_{+0}}f_{_{i+1,j}}+a_{_{-0}}f_{_{i-1,j}}
+a_{_{0+}}f_{_{i,j+1}}+a_{_{0-}}f_{_{i,j-1}}
\phantom{MMMMMM}\right.
\\
\phantom{MMMMMMMMMMMMMMMM}\left.
{}-\left(a_{_{+0}}+a_{_{-0}}+a_{_{0+}}+a_{_{0-}}\right)f_{_{i,j}}\right]
\end{array}
$$
and
$$
\begin{array}{l}
\L_{\times}^a[f]\equiv\dfrac{1}{2h^2}\left[
a_{_{++}}f_{_{i+1,j+1}}+a_{_{-+}}f_{_{i-1,j+1}}
+a_{_{+-}}f_{_{i+1,j-1}}+a_{_{--}}f_{_{i-1,j-1}}
\phantom{MMMMMM}\right.
\\
\phantom{MMMMMMMMMMMMMMMM}\left.
{}-\left(a_{_{++}}\!+a_{_{-+}}\!+a_{_{+-}}\!+a_{_{--}}\right)f_{_{i,j}}\right]\,,
\end{array}
$$
where $a_{_{\pm0}}$, $a_{_{0\pm}}$, and $a_{_{\pm\pm}}$  denote the values of $a(x,y)$
at $\left((i\pm 1/2)h,jh\right)$, $\left(ih,(j\pm 1/2)h\right)$, and
$\left((i\pm 1/2)h,(j\pm 1/2)h\right)$, respectively.

Assume $h\ll1$. Since the same weights are used in the discrete directional
mean value representations (\ref{dmv4lap}) and (\ref{qLap}) we can expect that
the optimal rotation-invariant  stencil for quasi-Laplacian $\L^a$
is obtained when $\alpha=2/3$ and $\beta=1/3$, or equivalently, when $w=4$.
This indeed turns out to be the case in the lower order terms:

\begin{proposition}
The discrete approximation $\alpha \L^a_{+}+\beta \L^a_{\times}$
of the quasi-Laplacian $\L^a$  is is asymptotically rotation-equivariant
for $\alpha=2/3$ and $\beta=1/3$
up to order $O\left(h^3\right)$ as $h\to 0$.
\end{proposition}

\noindent
{\sc Proof.}
A straightforward calculation shows that
\begin{equation}\label{aTaylor}
\begin{array}{l}
L^a[f]=\nabla\cdot[a(x,y)\nabla f]
\phantom{MMMMMM}
\vspace*{6pt}
\\
\phantom{MM}
{}+\dfrac{h^2}{12}\left[D_{_{0,4}}(a,f)+2D_{_{1,3}}(a,f)+\dfrac32D_{_{2,2}}(a,f)
+\dfrac12D_{_{3,1}}(a,f)\right]+O\left(h^4\right)
\end{array}
\end{equation}
with
$$
\begin{array}{rcl}
D_{_{0,4}}(a,f)&=&a(x,y)\left(
\dfrac{\partial^4f}{\partial x^4}+2\dfrac{\partial^4f}{\partial x^2\partial y^2}
+\dfrac{\partial^4f}{\partial y^4}
\right),
\vspace*{6pt}\\
D_{_{1,3}}(a,f)&=&
\dfrac{\partial a}{\partial x}\dfrac{\partial^3f}{\partial x^3}
+\dfrac{\partial a}{\partial y}\dfrac{\partial^3f}{\partial x^2\partial y}
+\dfrac{\partial a}{\partial x}\dfrac{\partial^3f}{\partial x\partial y^2}
+\dfrac{\partial a}{\partial y}\dfrac{\partial^3f}{\partial y^3},
\vspace*{6pt}\\
D_{_{2,2}}(a,f)&=&
\dfrac{\partial^2 a}{\partial x^2}\dfrac{\partial^2f}{\partial x^2}
+\dfrac{\partial^2 a}{\partial y^2}\dfrac{\partial^2f}{\partial y^2}
+\dfrac43\dfrac{\partial^2 a}{\partial x\partial y}\dfrac{\partial^2f}{\partial x\partial y}
+\dfrac13\dfrac{\partial^2 a}{\partial x^2}\dfrac{\partial^2f}{\partial y^2}
+\dfrac13\dfrac{\partial^2 a}{\partial y^2}\dfrac{\partial^2f}{\partial x^2},
\vspace*{6pt}\\
D_{_{3,1}}(a,f)&=&
\dfrac{\partial^3a}{\partial x^3}\dfrac{\partial f}{\partial x}
+\dfrac{\partial^3a}{\partial x^2\partial y}\dfrac{\partial f}{\partial y}
+\dfrac{\partial^3a}{\partial x\partial y^2}\dfrac{\partial f}{\partial x}
+\dfrac{\partial^3a}{\partial y^3}\dfrac{\partial f}{\partial y}.
\end{array}
$$
It is interesting that expansion (\ref{aTaylor}) does not contain $h^3$-terms.

Let us now demonstrate that these differential quantities
are rotationally invariant. Obviously $D_{_{0,4}}(a,f)$
equals $a(x,y)\Delta^2f$ which is rotationally invariant.
Let us consider
$$
R_{_{1,3}}(a,f,\varphi)=
\dfrac{\partial a}{\partial e_\varphi}\cdot\dfrac{\partial^3 f}{\partial e_\varphi^3},
\quad R_{_{2,2}}(a,f,\varphi)=
\dfrac{\partial^2 a}{\partial e_\varphi^2}\cdot\dfrac{\partial^2 f}{\partial e_\varphi^2},
\quad R_{_{3,1}}(a,f,\varphi)=
\dfrac{\partial^3 a}{\partial e_\varphi^3}\cdot\dfrac{\partial f}{\partial e_\varphi}.
$$
Note that $R_{_{3,1}}(a,f,\varphi)=R_{_{1,3}}(f,a,\varphi)$ and
$D_{_{1,3}}(a,f)=D_{_{3,1}}(f,a)$.
Direct calculations show that
$$
\int_0^{2\pi}R_{_{1,3}}(a,f,\varphi)\,d\varphi=\dfrac{3\pi}{4}D_{_{1,3}}(a,f)
\mbox{ and }
\int_0^{2\pi}R_{_{2,2}}(a,f,\varphi)\,d\varphi=\dfrac{3\pi}{4}D_{_{2,2}}(a,f).
$$
This completes our proof that the $h^2$-term in (\ref{aTaylor}) is rotation invariant. QED

\begin{figure}[htbp]
\centering
\includegraphics[width=0.29\textwidth]{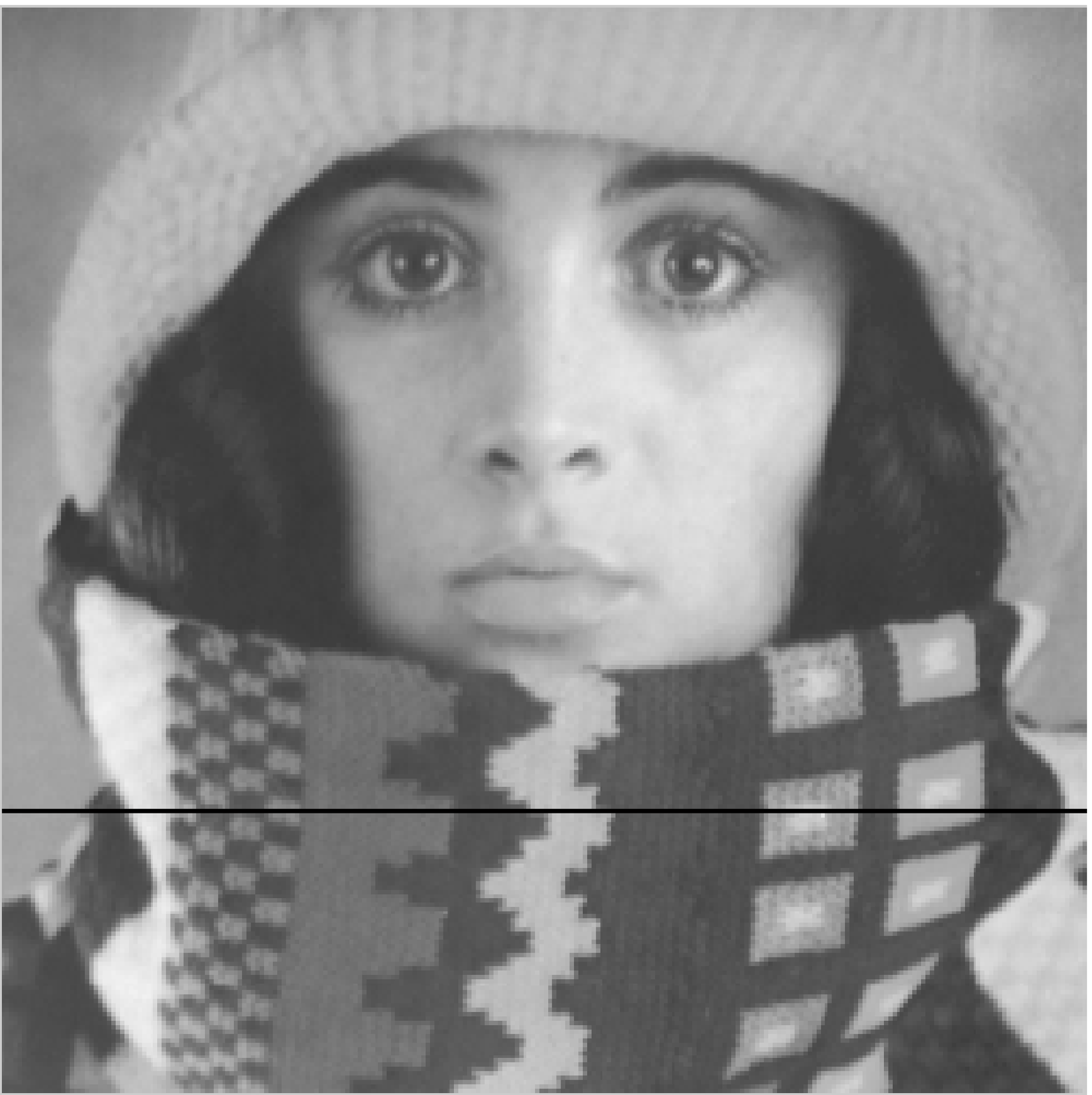}
\includegraphics[width=0.37\textwidth]{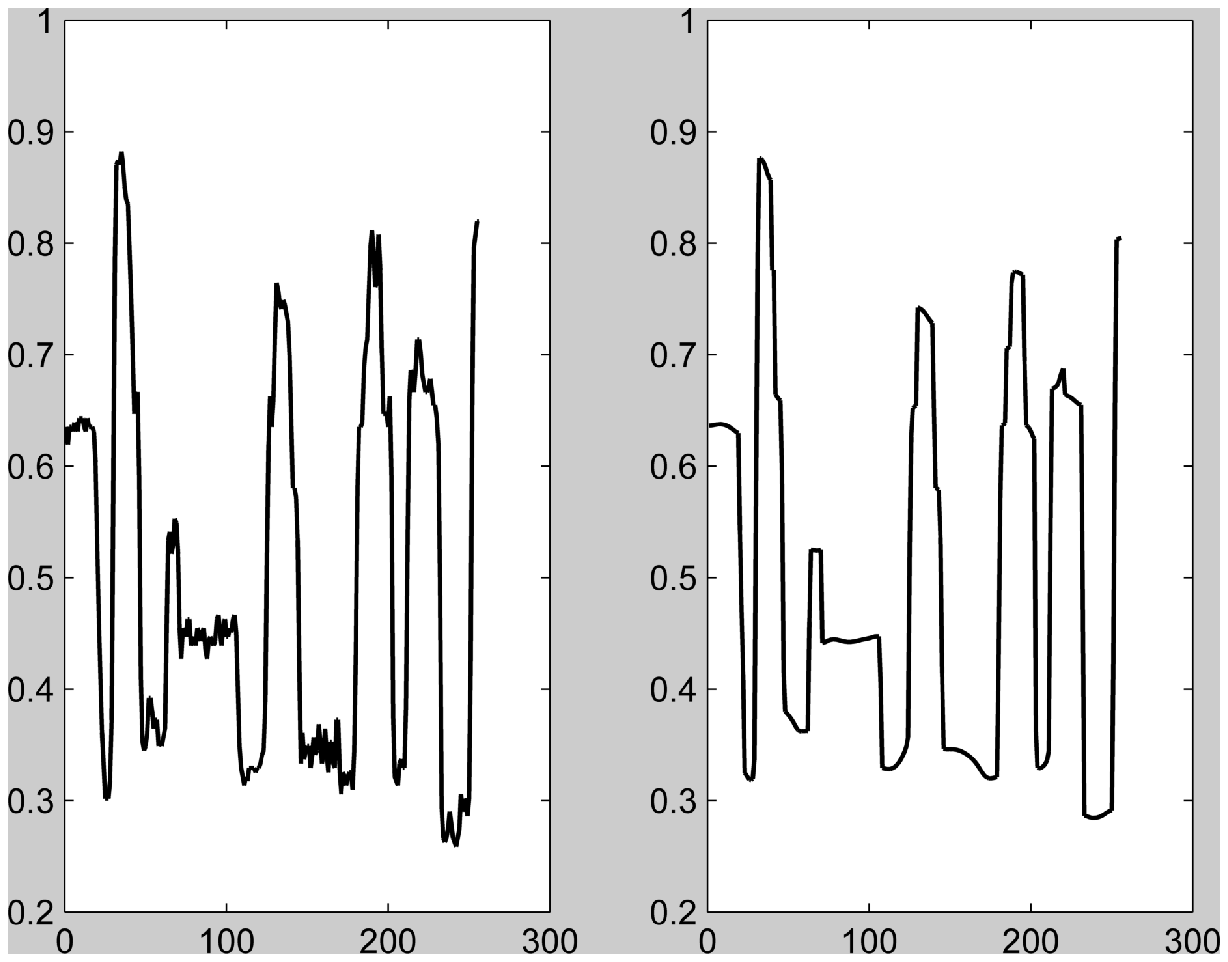}
\includegraphics[width=0.29\textwidth]{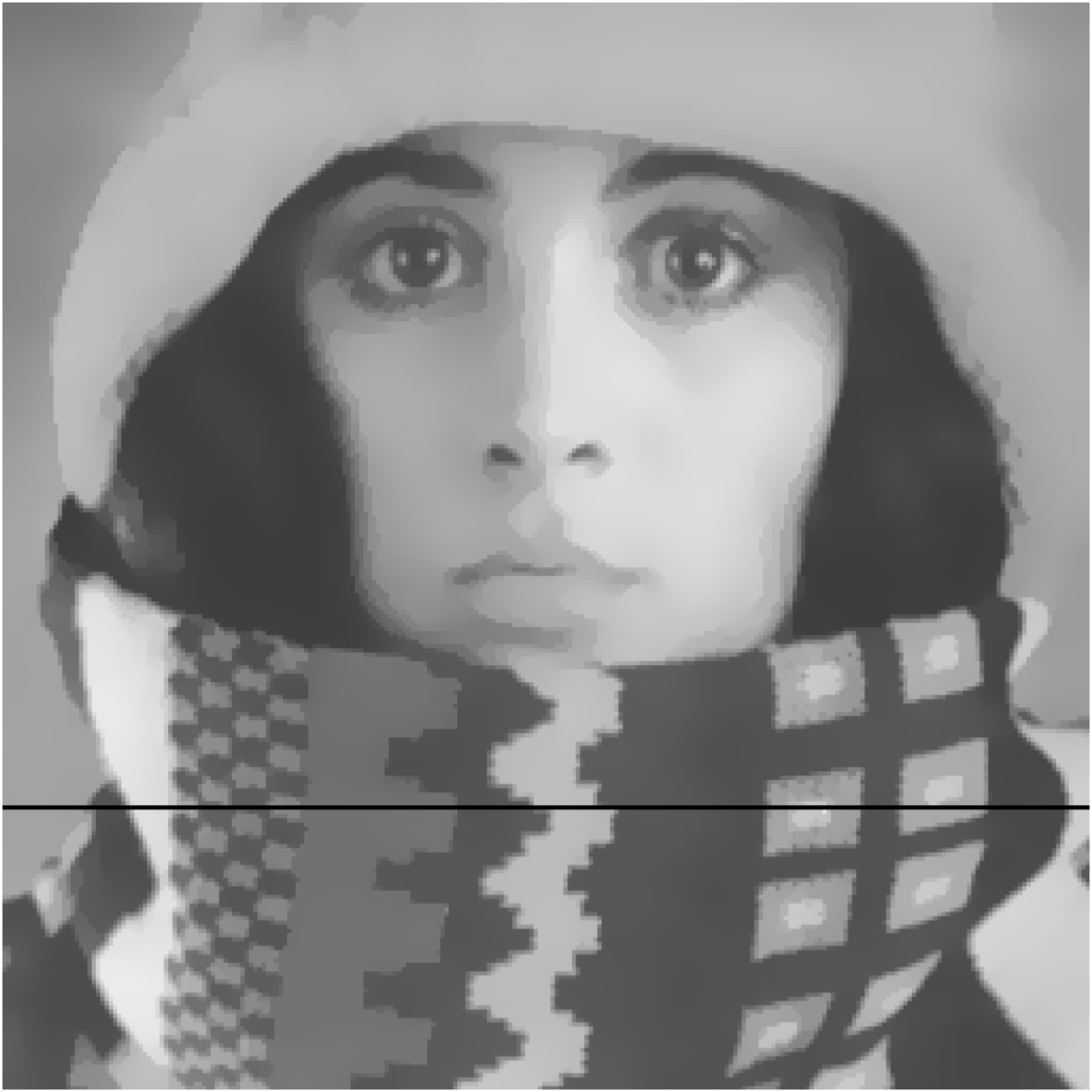}
\caption{Demonstrating the power of nonlinear diffusion for image filtering.
Left: original image ``trui". Right: after filtering with (\ref{nd}) and (\ref{nd-pm}).
Middle: a comparison of one-dimensional image slices of the original and filtered
images.}
\label{trui}
\end{figure}

\paragraph*{\bf Application to nonlinear diffusion image filtering}
Accurate numerical implementations of nonlinear diffusion processes governed by
\begin{equation}\label{nd}
\dfrac{\partial u(x,y,t)}{\partial t}
=\div\left(a(x,y,u,\nabla u)\nabla u\right),\quad u(x,y,0)=u_0(x,y),
\end{equation}
subject to appropriate boundary conditions, are important for a number of disciplines
including computational physics, magnetohydrodynamics, financial mathematics, and image processing.
In particular, in image processing, adaptive image smoothing is often carried out
by (\ref{nd}) with
\begin{equation}\label{nd-pm}
a(x,y,u,\nabla u)=\exp(-|\nabla u|/\lambda),
\end{equation}
as suggested by Perona and Malik in their seminal paper \cite{Perona-Malik_pami90}.
Figure~\ref{trui} demonstrates how our numerical implementation of
(\ref{nd}) and  (\ref{nd-pm}) performs an adaptive image smoothing and removes
texture from an image.

In image processing, the system (\ref{nd}),(\ref{nd-pm}) is typically solved
by finite differences methods with explicit or implicit Euler schemes and
several dozens of iterations are usually required in order to achieve desired
image filtering effects. For a number of applications including image segmentation
and feature extraction, accurate estimation of the gradient directions in an image
filtered by (\ref{nd}) and (\ref{nd-pm}) is required. Unfortunately standard finite
difference schemes have poor rotation-invariance properties. If such schemes are used
repeatedly, the directional error is accumulated. See
\cite{Weickert_ecmi94,Weickert_jvcir02} for examples of finite difference schemes
with improved rotation invariance properties.

\begin{figure}[h]
\centering
\includegraphics[width=0.32\textwidth]{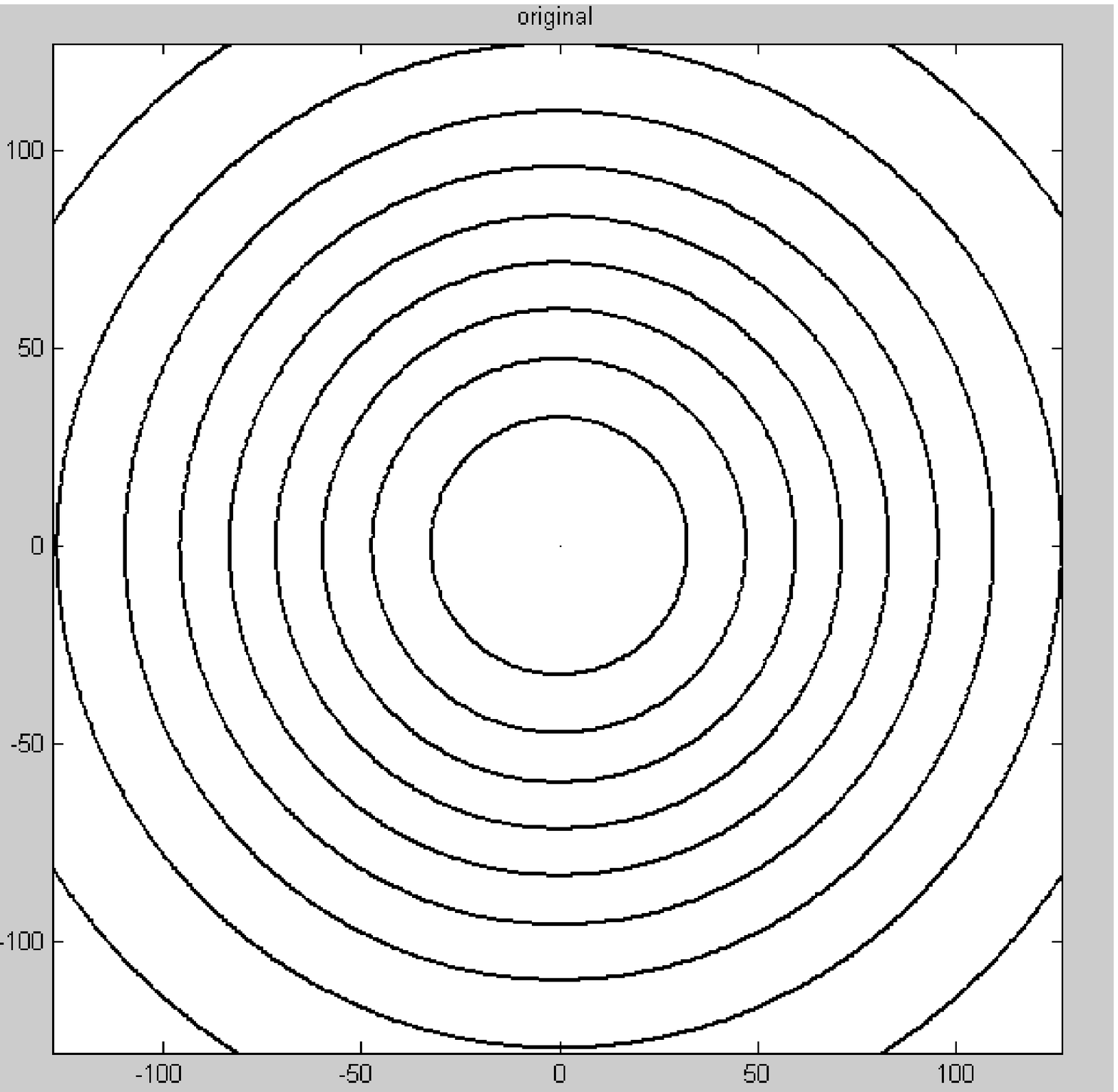}
\includegraphics[width=0.32\textwidth]{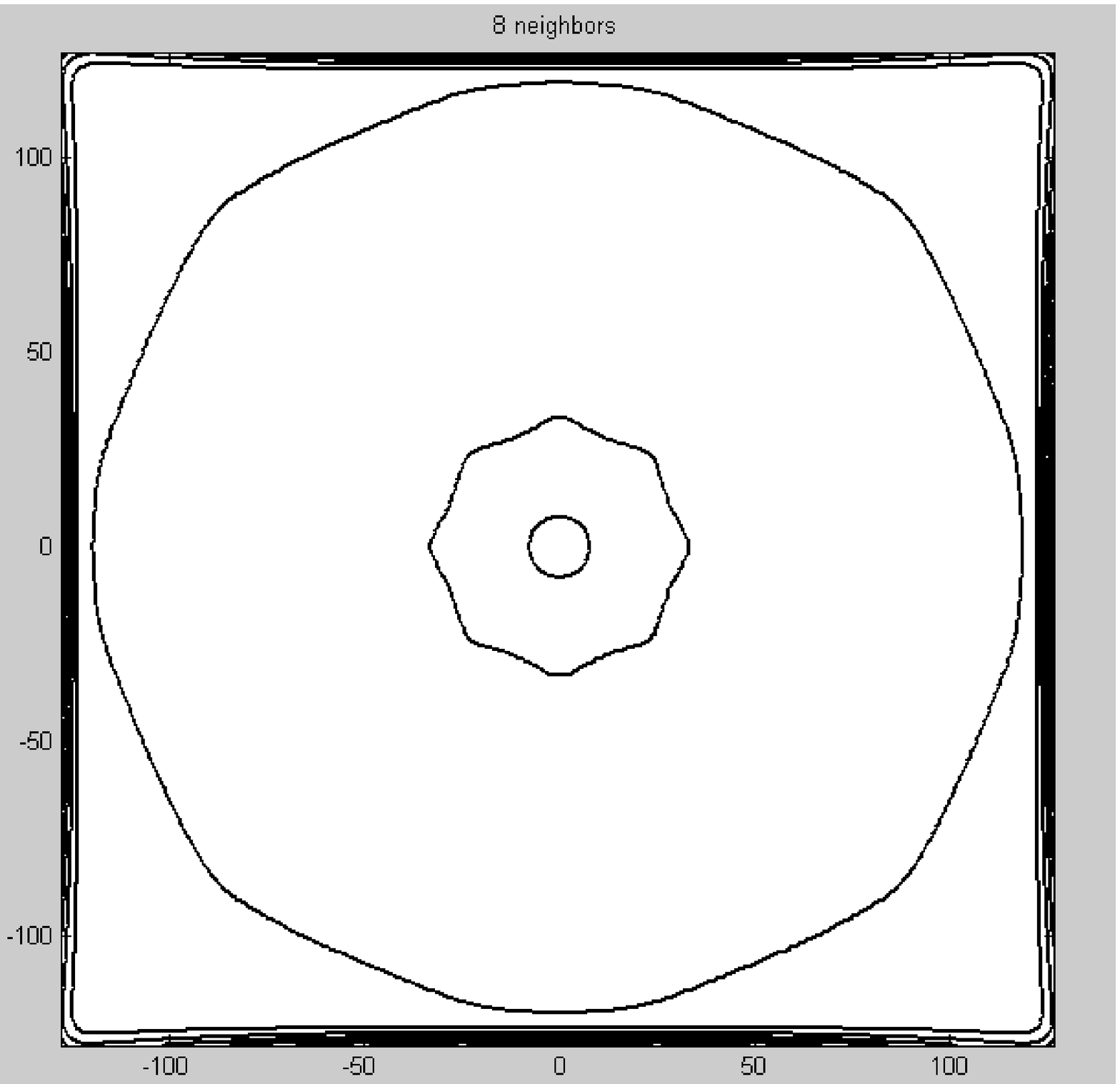}
\includegraphics[width=0.32\textwidth]{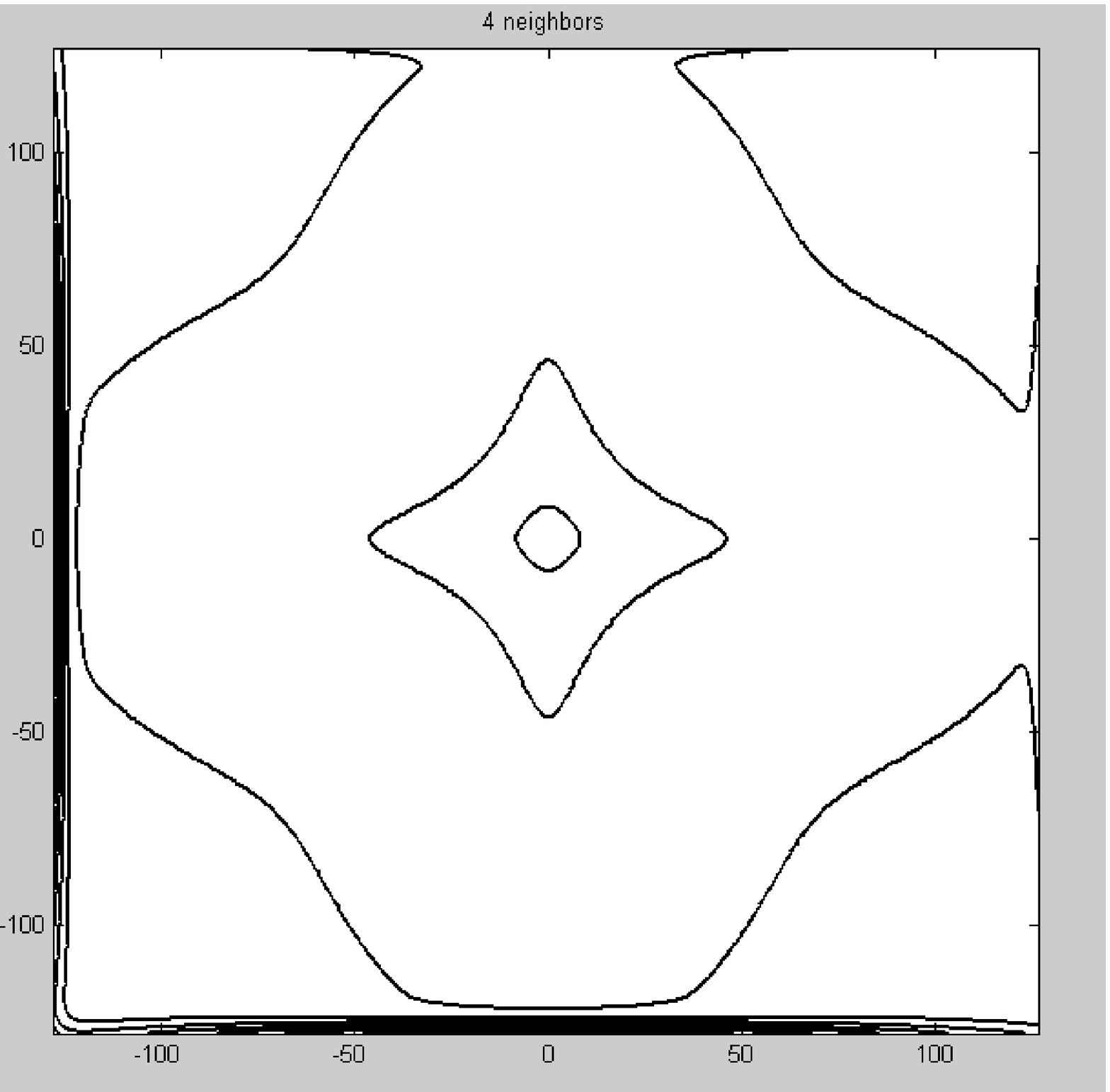}
\caption{Left: the contour lines of Gaussian bump function (\ref{gaussian})
are perfect circles. Middle: contours of the Gaussian filtered using our
implementation of nonlinear diffusion (\ref{nd}), (\ref{nd-pm}). Right:
contours of the Gaussian filtered using code from \cite{Perona-Shiota-Malik_94}.
}
\label{nd-gauss}
\end{figure}

Our approach to numerical solving of (\ref{nd}), (\ref{nd-pm}) combines
the explicit Euler scheme with obvious generalizations of the optimal
rotation-invariant $3\times3$ stencils derived in previous sections for gradient,
Laplacian, and quasi-Laplacian. To demonstrate advantages of our approach we apply
our numerical implementation of (\ref{nd}), (\ref{nd-pm}) to (\ref{gaussian}).
Figure~\ref{nd-gauss} shows level sets of (\ref{gaussian}) before and after
applying nonlinear diffusion (\ref{nd}), (\ref{nd-pm}). We also compare our
implementation with that from from \cite{Perona-Shiota-Malik_94}, where five-point
stencils are used for numerical solving (\ref{nd}) and (\ref{nd-pm}).


\paragraph*{\bf Frequency response analysis}
Now we take a brief look at stencils (\ref{wgrad}) and (\ref{wLap})
from a frequency point of view. This standpoint is important
when dealing with phenomena possessing a range
of space and time scales. Turbulent fluid flows, electromagnetic and
acoustic wave propagation, multiresolution analysis of financial
data, and multiscale image processing are among common examples.

For the sake of simplicity assume that our square grid
has a unit step-size: $h=1$ (this can be easily achieved by rescaling).
One can see that (\ref{wgrad}) consists of the simplest central difference operator combined with a smoothing kernel applied to the orthogonal direction:
\begin{equation}\label{Dx_split}
D_x=\frac{1}{2}\left[\begin{array}{rcl}
-1 &  0 & 1
\end{array}\right]
\cdot
\frac{1}{w+2}\left[\begin{array}{rcl}
1 & w & 1
\end{array}\right]^T.
\end{equation}
The eigenvalue (also called the frequency
response) of (\ref{wgrad}) corresponding
to the eigenfunction $\exp(i(\omega_1x+\omega_2y)$ is given by
$$
i\sin\omega_1\cdot\frac{w+2}{w+2\cos\omega_2},\quad
-\pi<\omega_k<\pi,\,\,\,k=1,2,
$$
where the first term of the product is the frequency response for
the central difference and the second one corresponds to the smoothing
kernel in (\ref{Dx_split}).

It is well known that the frequency response $i\sin\omega$ delivers
a satisfactory approximation of $i\,\omega$, the frequency response for
the ideal derivative, only for sufficiently small frequencies $\omega$
(see, for example, \cite[Section\,6.4]{Hamming_df}). Now it is clear
how (\ref{wgrad}) improves the central difference: smoothing due to the use
the central difference operator instead of the true $x$-derivative is
compensated by adding a certain amount of smoothing in the $y$-direction.
Thus (\ref{wgrad}) and its $y$-direction counterpart do a better job in estimating the gradient direction than in estimating the gradient magnitude.

If the goal is to achieve an accurate estimation of both the gradient direction and magnitude, we can combine (\ref{wgrad}) and (\ref{wLap}) as follows.
Let
$$
\delta = 
\left[\begin{array}{rcl}
0 & 0 & 0
\\
0 & 1 & 0
\\
0 & 0 & 0
\end{array}\right]
$$
be the $3\times3$ identity kernel. Note that
$$
\delta+\frac{1}{w+2}\,\wLap=
\frac{1}{(w+2)^2}
\left[\begin{array}{rcl}
1 & w & 1
\\
w & w^2 & w
\\
1 & w & 1
\end{array}\right]
=\frac{1}{w+2}\left[\begin{array}{rcl}
1 & w & 1
\end{array}\right]\cdot
\frac{1}{w+2}\left[\begin{array}{c}
1 \\ w \\ 1
\end{array}\right]\,,
$$
which can be considered as simultaneous smoothing (averaging) with respect to both the coordinate directions. Thus it is natural to use
\begin{equation}\label{sharp_wgrad}
\left(\delta+\frac{1}{w+2}\,\wLap\right)^{\!-1}\!\!D_x
\end{equation}
which combines (\ref{wgrad}) with a Laplacian-based sharpening.
The frequency response function corresponding to (\ref{sharp_wgrad})
is given by
\begin{equation}\label{frf_wgrad}
H(\omega)=i\sin\omega\cdot\frac{w+2}{w+2\cos\omega}\,,
\end{equation}
where parameter $w$ should be chosen in such a way that (\ref{frf_wgrad})
delivers a good approximation of $i\,\omega$, the frequency response for
the true derivative. For example, the Taylor series expansion of
(\ref{frf_wgrad}) with respect to  $\omega$ at $\omega=0$ shows that the best
approximation of the exact $x$-derivative for $\omega\ll1$ is
achieved when $w=4$ which corresponds to a Pad{\'e} approximation.
However this choice of parameter $w$ is not
optimal if we deal with a wider range of frequencies.

Note that (\ref{sharp_wgrad}) is equivalent to an
implicit finite difference scheme
\begin{equation}\label{ifd1}
\frac{1}{w+2}
\left(f^\prime_{i-1,j}+wf^\prime_{i,j}+f^\prime_{i+1,j}\right)
=\frac12\left(f_{i+1,j}-f_{i-1,j}\right).
\end{equation}
Implicit finite differences are widely used for accurate numerical simulations
of physical problems involving linear and non-linear wave propagation phenomena,
see, for example,
\cite{Colonius-Lele_pas04}, \cite[Section\,5.8]{Petrila-Trif_cfd05}.
An implicit finite difference scheme is usually characterized by
its frequency-resolving efficiency, the range of frequencies $\omega$
over which a satisfactory approximation of the corresponding differential operator is achieved.

\begin{figure}[htbp]
\centering
\includegraphics[width=0.6\textwidth]{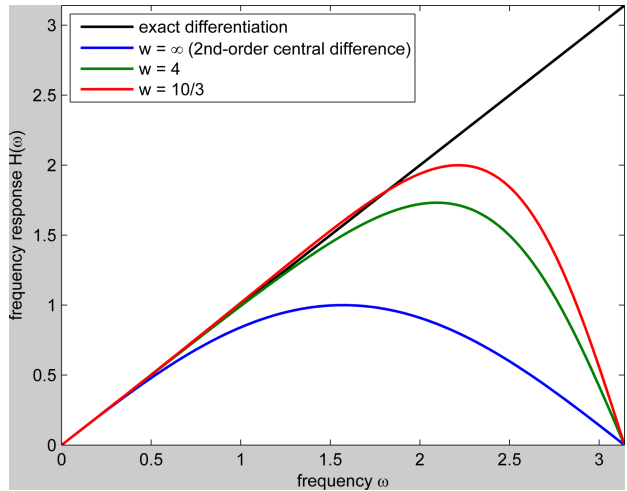}
\caption{
Visual frequency-response analysis of (\ref{sharp_wgrad})
for various values of parameter $w$.
Notice good frequency-resolving efficiency of (\ref{frf_wgrad})
with $w=10/3$.
}
\label{implicit}
\end{figure}

Motivated by \cite{Scharr_dagm97,Jahne_sppr99} where a frequency-based
analysis of simple $3\times3$ stencils for first-order derivatives was
performed, we found out that (\ref{frf_wgrad}, (\ref{ifd1})
with $w=10/3$ deliver a good frequency resolving
efficiency for a sufficiently large frequency range, as seen in Figure~\ref{implicit}. Of course,
the problem of finding an optimal $w$ is application-dependent and
the optimal value of $w$ depends on optimization criteria used and
a range of frequencies where the optimization is sought.

\medskip

The approach via the Minkowski theorem and Veronese map developed in Section \ref{sec-t}
is suggestive for an extension of the above results to 3D and to grids in any dimensions. We leave this for future publication.

\medskip

{\bf Acknowledgments}.  B.\,K. is grateful to the Ecole Polytechnique in Paris and the Max
Planck Institute in Bonn for their hospitality during completion of
this paper.  He was partially supported by an NSERC research grant.
S.\,T. was partially supported by the Simons Foundation grant No 209361 and by the NSF grant DMS-1105442.

\bibliographystyle{elsarticle-harv}
\bibliography{derivatives}

\begin{thebibliography}{40}
\expandafter\ifx\csname natexlab\endcsname\relax\def\natexlab#1{#1}\fi
\expandafter\ifx\csname url\endcsname\relax
  \def\url#1{\texttt{#1}}\fi
\expandafter\ifx\csname urlprefix\endcsname\relax\def\urlprefix{URL }\fi

\bibitem[{Arnold and Khesin(1998)}]{Arnold-Khesin98}
Arnold, V., Khesin, B., 1998. Topological methods in hydrodynamics. Springer.

\bibitem[{Aubert and Kornprobst(2006)}]{Aubert-Kornprobst_mpip2e06}
Aubert, G., Kornprobst, P., 2006. Mathematical Problems in Image Processing,
  2nd Edition. Springer.

\bibitem[{Bannai and Bannai(2009)}]{Bannai}
Bannai, E., Bannai, E., 2009. A survey on spherical designs and algebraic
  combinatorics on spheres. European Journal of Combinatorics 30~(6),
  1392--1425.

\bibitem[{Belyaev(2006)}]{Belyaev_sgp06}
Belyaev, A., June 2006. On transfinite barycentric coordinates. In: Proceedings
  of Fourth Eurographics Symposium on Geometry Processing. Sardinia, Italy, pp.
  89--99.

\bibitem[{Bickley(1948)}]{Bickley_qjmam48}
Bickley, W.~G., 1948. Finite difference formulae for the square lattice. Quart.
  J. Mech. Appl. Math. 1, 35--42.

\bibitem[{Birkhoff and Lynch(1984)}]{Birkhoff-Lynch_nsep84}
Birkhoff, G., Lynch, R.~E., 1984. Numerical Solution of Elliptic Problems.
  SIAM.

\bibitem[{Bradford and Westhead(2005)}]{Bradford-Westhead_bio05}
Bradford, J.~R., Westhead, D.~R., 2005. Improved prediction of protein-protein
  binding sites using a support vector machines approach. Bioinformatics
  21~(8), 1487--1494.

\bibitem[{Childress(1992)}]{Childress}
Childress, S., 1992. Fast dynamo theory. In: Topological aspects of the
  dynamics of fluids. Kluwer Academic Publishers, Dordrecht, pp. 111--147.

\bibitem[{Colonius and Lele(2004)}]{Colonius-Lele_pas04}
Colonius, T., Lele, S.~K., August 2004. Computational aeroacoustics: progress
  on nonlinear problems of sound generation. Progress in Aerospace Sciences
  40~(6), 345--416.

\bibitem[{Delsarte et~al.(1977)Delsarte, Goethals, and Seidel}]{DGS}
Delsarte, P., Goethals, J.~M., Seidel, J.~J., 1977. Spherical codes and
  designs. Geometriae Dedicata 6, 363--388.

\bibitem[{Friedman and Littman(1962)}]{Friedman-Littman_62}
Friedman, A., Littman, W., 1962. Functions satisfying the mean value property.
  Transactions of the American Mathematical Society 102~(1), 167--180.

\bibitem[{Gonzalez and Woods(2008)}]{Gonzalez-Woods_dip3e08}
Gonzalez, R.~C., Woods, R.~E., 2008. Digital Image Processing, 3rd Edition.
  Pearson Prentice Hall.

\bibitem[{Gordon and Wixom(1974)}]{Gordon-Wixom_sinum74}
Gordon, W., Wixom, J., 1974. Pseudo-harmonic interpolation on convex domains.
  SIAM J. Numer. Anal. 11~(5), 909--933.

\bibitem[{Hamming(1998)}]{Hamming_df}
Hamming, R.~W., 1998. Digital Filters, 3rd Edition. Dover.

\bibitem[{Harris(1992)}]{Harris_ag}
Harris, J., 1992. Algebraic Geometry: A First Course. Springer.

\bibitem[{Hong(1982)}]{Hong-1982}
Hong, Y., 1982. On spherical t-designs in $r^2$. European Journal of
  Combinatorics 3, 255--258.

\bibitem[{J{\"a}hne et~al.(1999)J{\"a}hne, Scharr, and
  K{\"o}rkel}]{Jahne_sppr99}
J{\"a}hne, B., Scharr, H., K{\"o}rkel, S., 1999. Principles of filter design.
  In: Handbook of Computer Vision and Applications, Vol.\,2, Signal Processing
  and Applications. Academic Press, pp. 125--151.

\bibitem[{Kamgar-Parsi et~al.(1999)Kamgar-Parsi, Kamgar-Parsi, and
  Rosenfeld}]{Kamgar-Parsi_tip99}
Kamgar-Parsi, B., Kamgar-Parsi, B., Rosenfeld, A., 1999. Optimally isotropic
  {L}aplacian operator. IEEE Transactions on Image Processing 8~(10),
  1467--1472.

\bibitem[{Kantorovich and Krylov(1956)}]{Kantorovich-Krylov_amha56}
Kantorovich, L.~V., Krylov, V.~I., 1956. Approximate Methods of Higher
  Analysis. Noordhoff-Interscience.

\bibitem[{Koenderink and Van~Doorn(1992)}]{Koenderink-VanDoorn_92}
Koenderink, J.~J., Van~Doorn, A.~J., 1992. Surface shape and curvature scales.
  Image and Vision Computing 10, 557--565.

\bibitem[{Langer et~al.(2007)Langer, Belyaev, and Seidel}]{Langer_cagd07}
Langer, T., Belyaev, A., Seidel, H.-P., 2007. Exact and interpolatory
  quadratures for curvature tensor estimation. Computer Aided Geometric Design
  24~(8-9), 443--463.

\bibitem[{M{\'e}rigot et~al.(2009)M{\'e}rigot, Ovsjanikov, and
  Guibas}]{Meigot_spm09}
M{\'e}rigot, Q., Ovsjanikov, M., Guibas, L., 2009. Robust {V}oronoi-based
  curvature and feature estimation. In: Proc. of SIAM/ACM Joint Conference on
  Geometric and Physical Modeling. pp. 1--12.

\bibitem[{Natterer(2001)}]{Natterer_siam01}
Natterer, F., 2001. The Mathematics of Computerized Tomography. SIAM.

\bibitem[{Patra and Karttunen(2005)}]{Patra-Karttunen_nmpde05}
Patra, M., Karttunen, M., 2005. Stencils with isotropic discretisation error
  for differential operators. Numerical Methods for Partial Differential
  Equations 22~(4), 936--953.

\bibitem[{Perona and Malik(1990)}]{Perona-Malik_pami90}
Perona, P., Malik, J., 1990. Scale-space and edge detection using anisotropic
  diffusion. IEEE Transactions on Pattern Analysis and Machine Intelligence
  12~(7), 629--639.

\bibitem[{Perona et~al.(1994)Perona, Shiota, and
  Malik}]{Perona-Shiota-Malik_94}
Perona, P., Shiota, T., Malik, J., 1994. Anisotropic Diffusion. Kluwer Academic
  Press, pp. 73--92.

\bibitem[{Petrila and Trif(2005)}]{Petrila-Trif_cfd05}
Petrila, T., Trif, D., 2005. Basics of Fluid Mechanics and Introduction to
  Computational Fluid Dynamics. Springer.

\bibitem[{Pottmann et~al.(2009)Pottmann, Wallner, Huang, and
  Yang}]{pottmann-2009-iir}
Pottmann, H., Wallner, J., Huang, Q., Yang, Y.-L., 2009. Integral invariants
  for robust geometry processing. Computer Aided Geometric Design 26, 37--60.

\bibitem[{Pratt(2001)}]{Pratt_dip3e01}
Pratt, W.~K., 2001. Digital Image Processing, 3rd Edition. John Wiley \& Sons.

\bibitem[{Scharr et~al.(1997)Scharr, K{\"o}rkel, and J{\"a}hne}]{Scharr_dagm97}
Scharr, H., K{\"o}rkel, S., J{\"a}hne, B., 1997. Numerische
  isotropieoptimierung von fir-filtern mittels quergl{\"a}ttung. In: Proc.
  DAGM'97. pp. 367--374.

\bibitem[{Schneider(1993)}]{Schneider_93}
Schneider, R., 1993. Convex Bodies: The Brunn-Minkowski Theory. Cambridge
  University Press.

\bibitem[{Seymour and Zaslavsky(1984)}]{SZ}
Seymour, P.~D., Zaslavsky, T., 1984. Averaging sets: {A} generalization of mean
  values and spherical designs. Advances in Math. 52, 213–224.

\bibitem[{Stroud(1971)}]{Stroud_book71}
Stroud, A.~H., 1971. Approximate Calculation of Multiple Integrals.
  Prentice-Hall.

\bibitem[{Taubin(1995)}]{Taubin_iccv95}
Taubin, G., 1995. Estimating the tensor of curvature of a surface from a
  polyhedral approximation. In: Proc. ICCV'95. pp. 902--907.

\bibitem[{Valenti et~al.(2009)Valenti, Sebe, and Gevers}]{ValentiICCV2009}
Valenti, R., Sebe, N., Gevers, T., 2009. Image saliency by isocentric
  curvedness and color. In: IEEE International Conference on Computer Vision.

\bibitem[{Wardetzky et~al.(2007)Wardetzky, Mathur, K{\"a}lberer, and
  E.}]{Wardetzky_sgp07}
Wardetzky, M., Mathur, S., K{\"a}lberer, F., E., G., July 2007. Discrete
  {L}aplace operators: {N}o free lunch. In: Proceedings of Fifth Eurographics
  Symposium on Geometry Processing. Barcelona, Spain, pp. 33--37.

\bibitem[{Weickert(1994)}]{Weickert_ecmi94}
Weickert, J., 1994. Anisotropic diffusion filters for image processing based
  quality control. In: Proc. Seventh European Conf. on Mathematics in Industry.
  Teubner-Verlag, pp. 355--362.

\bibitem[{Weickert(1998)}]{Weickert_book98}
Weickert, J., 1998. Anisotropic Diffusion in Image Processing. Teubner-Verlag,
  Stuttgart.

\bibitem[{Weickert and Scharr(2002)}]{Weickert_jvcir02}
Weickert, J., Scharr, H., 2002. A scheme for coherence-enhancing diffusion
  filtering with optimized rotation invariance. Journal of Visual Communication
  and Image Representation 13, 103--118.

\bibitem[{Zhong et~al.(2009)Zhong, Su, Yeo, Tan, Ghista, and Kassab}]{Zhong_09}
Zhong, L., Su, Y., Yeo, S.-Y., Tan, R.-S., Ghista, D.~N., Kassab, G., 2009.
  Left ventricular regional wall curvedness and wall stress in patients with
  ischemic dilated cardiomyopathy. American Journal of Physiology - Heart and
  Circulatory Physiology 296, H573--H584.

\end{thebibliography}
\end{document}